\documentclass[10pt]{amsart}
\usepackage{amsfonts}
\usepackage{amsmath}
\usepackage{amsthm}
\usepackage{amssymb}
\usepackage{latexsym}
\newtheorem{cor}{Corollary}[section]
\newtheorem{lem}{Lemma}[section]

\newtheorem{prop}{Proposition}[section]

\theoremstyle{definition}
\newtheorem{defn}{Definition}[section]
\theoremstyle{definition}
\newtheorem{thm}{Theorem}

\newtheorem*{rem}{Remark}

\newenvironment{pf}{\proof}{\endproof}

\theoremstyle{remark}

\numberwithin{equation}{section}
\begin{document}

\newcommand{\thmref}[1]{Theorem~\ref{#1}}
\newcommand{\secref}[1]{Sect.~\ref{#1}}
\newcommand{\lemref}[1]{Lemma~\ref{#1}}
\newcommand{\propref}[1]{Proposition~\ref{#1}}
\newcommand{\corref}[1]{Corollary~\ref{#1}}
\newcommand{\remref}[1]{Remark~\ref{#1}}
\newcommand{\er}[1]{(\ref{#1})}
\newcommand{\nc}{\newcommand}
\newcommand{\rnc}{\renewcommand}
\nc{\cal}{\mathcal}
\nc{\goth}{\mathfrak}
\rnc{\bold}{\mathbf}
\renewcommand{\frak}{\mathfrak}
\renewcommand{\Bbb}{\mathbb}

\nc{\Cal}{\cal}
\nc{\Xp}[1]{X^+(#1)}
\nc{\Xm}[1]{X^-(#1)}
\nc{\on}{\operatorname}
\nc{\ch}{\mbox{ch}}
\nc{\Z}{{\bold Z}}
\nc{\J}{{\cal J}}
\nc{\C}{{\bold C}}
\nc{\Q}{{\bold Q}}
\renewcommand{\P}{{\cal P}}
\nc{\N}{{\Bbb N}}
\nc\beq{\begin{equation}}
\nc\enq{\end{equation}}
\nc\lan{\langle}
\nc\ran{\rangle}
\nc\bsl{\backslash}
\nc\mto{\mapsto}
\nc\lra{\leftrightarrow}
\nc\hra{\hookrightarrow}
\nc\sm{\smallmatrix}
\nc\esm{\endsmallmatrix}
\nc\sub{\subset}
\nc\ti{\tilde}
\nc\nl{\newline}
\nc\fra{\frac}
\nc\und{\underline}
\nc\ov{\overline}
\nc\ot{\otimes}
\nc\bbq{\bar{\bq}_l}
\nc\bcc{\thickfracwithdelims[]\thickness0}
\nc\ad{\text{\rm ad}}
\nc\Ad{\text{\rm Ad}}
\nc\Hom{\text{\rm Hom}}
\nc\End{\text{\rm End}}
\nc\Ind{\text{\rm Ind}}
\nc\Res{\text{\rm Res}}
\nc\Ker{\text{\rm Ker}}
\rnc\Im{\text{Im}}
\nc\sgn{\text{\rm sgn}}
\nc\tr{\text{\rm tr}}
\nc\Tr{\text{\rm Tr}}
\nc\supp{\text{\rm supp}}
\nc\card{\text{\rm card}}
\nc\bst{{}^\bigstar\!}
\nc\he{\heartsuit}
\nc\clu{\clubsuit}
\nc\spa{\spadesuit}
\nc\di{\diamond}

\nc\al{\alpha}
\nc\bet{\beta}
\nc\ga{\gamma}
\nc\de{\delta}
\nc\ep{\epsilon}
\nc\io{\iota}
\nc\om{\omega}
\nc\si{\sigma}
\rnc\th{\theta}
\nc\ka{\kappa}
\nc\la{\lambda}
\nc\ze{\zeta}

\nc\vp{\varpi}
\nc\vt{\vartheta}
\nc\vr{\varrho}

\nc\Ga{\Gamma}
\nc\De{\Delta}
\nc\Om{\Omega}
\nc\Si{\Sigma}
\nc\Th{\Theta}
\nc\La{\Lambda}

\nc\boa{\bold a}
\nc\bob{\bold b}
\nc\boc{\bold c}
\nc\bod{\bold d}
\nc\boe{\bold e}
\nc\bof{\bold f}
\nc\bog{\bold g}
\nc\boh{\bold h}
\nc\boi{\bold i}
\nc\boj{\bold j}
\nc\bok{\bold k}
\nc\bol{\bold l}
\nc\bom{\bold m}
\nc\bon{\bold n}
\nc\boo{\bold o}
\nc\bop{\bold p}
\nc\boq{\bold q}
\nc\bor{\bold r}
\nc\bos{\bold s}
\nc\bou{\bold u}
\nc\bov{\bold v}
\nc\bow{\bold w}
\nc\boz{\bold z}

\nc\ba{\bold A}
\nc\bb{\bold B}
\nc\bc{\bold C}
\nc\bd{\bold D}
\nc\be{\bold E}
\nc\bg{\bold G}
\nc\bh{\bold H}
\nc\bi{\bold I}
\nc\bj{\bold J}
\nc\bk{\bold K}
\nc\bl{\bold L}
\nc\bm{\bold M}
\nc\bn{\bold N}
\nc\bo{\bold O}
\nc\bp{\bold P}
\nc\bq{\bold Q}
\nc\br{\bold R}
\nc\bs{\bold S}
\nc\bt{\bold T}
\nc\bu{\bold U}
\nc\bv{\bold V}
\nc\bw{\bold W}
\nc\bz{\bold Z}
\nc\bx{\bold X}

\nc\ca{\Cal A}
\nc\cb{\Cal B}
\nc\cc{\Cal C}
\nc\cd{\Cal D}
\nc\ce{\Cal E}
\nc\cf{\Cal F}
\nc\cg{\Cal G}
\rnc\ch{\Cal H}
\nc\ci{\Cal I}
\nc\cj{\Cal J}
\nc\ck{\Cal K}
\nc\cl{\Cal L}
\nc\cm{\Cal M}
\nc\cn{\Cal N}
\nc\co{\Cal O}
\nc\cp{\Cal P}
\nc\cq{\Cal Q}
\nc\car{\Cal R}
\nc\cs{\Cal S}
\nc\ct{\Cal T}
\nc\cu{\Cal U}
\nc\cv{\Cal V}
\nc\cz{\Cal Z}
\nc\cx{\Cal X}
\nc\cy{\Cal Y}

\nc\e[1]{E_{#1}}
\nc\ei[1]{E_{\delta - \alpha_{#1}}}
\nc\esi[1]{E_{s \delta - \alpha_{#1}}}
\nc\eri[1]{E_{r \delta - \alpha_{#1}}}
\nc\ed[2][]{E_{#1 \delta,#2}}
\nc\ekd[1]{E_{k \delta,#1}}
\nc\emd[1]{E_{m \delta,#1}}
\nc\erd[1]{E_{r \delta,#1}}

\nc\ef[1]{F_{#1}}
\nc\efi[1]{F_{\delta - \alpha_{#1}}}
\nc\efsi[1]{F_{s \delta - \alpha_{#1}}}
\nc\efri[1]{F_{r \delta - \alpha_{#1}}}
\nc\efd[2][]{F_{#1 \delta,#2}}
\nc\efkd[1]{F_{k \delta,#1}}
\nc\efmd[1]{F_{m \delta,#1}}
\nc\efrd[1]{F_{r \delta,#1}}

\nc\fa{\frak a}
\nc\fb{\frak b}
\nc\fc{\frak c}
\nc\fd{\frak d}
\nc\fe{\frak e}
\nc\ff{\frak f}
\nc\fg{\frak g}
\nc\fh{\frak h}
\nc\fj{\frak j}
\nc\fk{\frak k}
\nc\fl{\frak l}
\nc\fm{\frak m}
\nc\fn{\frak n}
\nc\fo{\frak o}
\nc\fp{\frak p}
\nc\fq{\frak q}
\nc\fr{\frak r}
\nc\fs{\frak s}
\nc\ft{\frak t}
\nc\fu{\frak u}
\nc\fv{\frak v}
\nc\fz{\frak z}
\nc\fx{\frak x}
\nc\fy{\frak y}

\nc\fA{\frak A}
\nc\fB{\frak B}
\nc\fC{\frak C}
\nc\fD{\frak D}
\nc\fE{\frak E}
\nc\fF{\frak F}
\nc\fG{\frak G}
\nc\fH{\frak H}
\nc\fJ{\frak J}
\nc\fK{\frak K}
\nc\fL{\frak L}
\nc\fM{\frak M}
\nc\fN{\frak N}
\nc\fO{\frak O}
\nc\fP{\frak P}
\nc\fQ{\frak Q}
\nc\fR{\frak R}
\nc\fS{\frak S}
\nc\fT{\frak T}
\nc\fU{\frak U}
\nc\fV{\frak V}
\nc\fZ{\frak Z}
\nc\fX{\frak X}
\nc\fY{\frak Y}
\nc\tfi{\ti{\Phi}}
\nc\bF{\bold F}
\rnc\bol{\bold 1}

\nc\ua{{_\ca\bold U}}

\nc\qinti[1]{[#1]_i}
\nc\q[1]{[#1]_q}
\nc\xpm[2]{E_{#2 \delta \pm \alpha_#1}}  
\nc\xmp[2]{E_{#2 \delta \mp \alpha_#1}}
\nc\xp[2]{E_{#2 \delta + \alpha_{#1}}}
\nc\xm[2]{E_{#2 \delta - \alpha_{#1}}}
\nc\hik{\ed{k}{i}}
\nc\hjl{\ed{l}{j}}
\nc\qcoeff[3]{\left[ \begin{smallmatrix} {#1}& \\ {#2}& \end{smallmatrix}
\negthickspace \right]_{#3}}
\nc\qi{q}
\nc\qj{q}

\nc\ufdm{{_\ca\bu}_{\rm fd}^{\le 0}}


\nc\isom{\cong} 

\nc{\pone}{{\Bbb C}{\Bbb P}^1}
\nc{\pa}{\partial}
\def\H{\cal H}
\def\L{\cal L}
\nc{\F}{{\cal F}}
\nc{\Sym}{{\goth S}}
\nc{\A}{{\cal A}}
\nc{\arr}{\rightarrow}
\nc{\larr}{\longrightarrow}

\nc{\ri}{\rangle}
\nc{\lef}{\langle}
\nc{\W}{{\cal W}}
\nc{\uqatwoatone}{{U_{q,1}}(\su)}
\nc{\uqtwo}{U_q(\goth{sl}_2)}
\nc{\dij}{\delta_{ij}}
\nc{\divei}{E_{\alpha_i}^{(n)}}
\nc{\divfi}{F_{\alpha_i}^{(n)}}
\nc{\Lzero}{\Lambda_0}
\nc{\Lone}{\Lambda_1}
\nc{\ve}{\varepsilon}
\nc{\phioneminusi}{\Phi^{(1-i,i)}}
\nc{\phioneminusistar}{\Phi^{* (1-i,i)}}
\nc{\phii}{\Phi^{(i,1-i)}}
\nc{\Li}{\Lambda_i}
\nc{\Loneminusi}{\Lambda_{1-i}}
\nc{\vtimesz}{v_\ve \otimes z^m}

\nc{\asltwo}{\widehat{\goth{sl}_2}}
\nc\ag{\widehat{\goth{g}}}  
\nc\teb{\tilde E_\boc}
\nc\tebp{\tilde E_{\boc'}}
\title{ An Algebraic characterization \\ of the affine canonical basis}
\author{Jonathan Beck} \
\address{Jonathan Beck, University of Chicago}
\author{Vyjayanthi Chari}
\address{Vyjayanthi Chari, University of California, Riverside}
\author{Andrew Pressley}
\address{Andrew Pressley, King's College, London}
\begin{abstract}
  The canonical basis for finite type quantized universal enveloping algebras was introduced
  in \cite{L1}. The principal technique is the explicit construction
  (via the braid group action) of a lattice $\cl$ over $\bz[q^{-1}]$.
  This allows the algebraic characterization of the canonical basis as a
  certain bar-invariant basis of $\cl$. Here we present a similar
  algebraic characterization of the affine canonical basis.  Our
  construction is complicated by the need to introduce basis elements to
  span the ``imaginary'' subalgebra which is fixed by the affine braid
  group.  Once the basis is found we construct a PBW-type basis whose
  $\bz[q^{-1}]$-span reduces to a ``crystal'' basis at $q=\infty,$
with the imaginary component given by the Schur functions.\end{abstract}
\pagestyle{plain} \thanks{Revised Version, June 1998.  V. C. is supported by a
  grant from NATO. A.P. is supported by the EPSRC (GR/L26216)}
\maketitle
\section{Introduction.}
The canonical basis of the  quantized universal enveloping algebra associated to a simple finite-dimensional Lie algebra was introduced by Lusztig in \cite{L1} via an
elementary algebraic definition.  The definition was characterized by
three main components: 1) the basis was integral, 2) it was
bar-invariant, and 3) it spanned a certain $\bz[q^{-1}]$-lattice $\cl$
with a specific image in the quotient $\cl/q^{-1}\cl.$ This algebraic
definition does not work for quantized
universal enveloping algebras of arbitrary Kac--Moody algebras. The difficulty in constructing a basis for $\cl$ arises from the need to define suitable analogues of imaginary root vectors.  The definition of the canonical basis for arbitrary
type was subsequently made using topological methods \cite{L5}. In his paper
\cite{Ka}, Kashiwara gave a suitable algebraic definition of the
lattices $\cl$ and $\cl/q^{-1}\cl$, making use of a remarkable symmetric bilinear form on the algebra, (introduced by Drinfeld) which led to an inductive construction of
the global crystal basis.  It was later shown in \cite{GL} that the two
concepts--the global crystal basis (algebraic) and the canonical basis
(topological)--coincide.

In this paper we synthesize the two aforementioned techniques and
construct a crystal basis for the quantized universal enveloping algebra
of (untwisted) affine type. Then we give an elementary algebraic
characterization of the canonical basis analogous to the
characterization given in the finite type  case  \cite{L1}. A remarkable feature of
our construction is that the part of the crystal basis corresponding to
the imaginary root spaces is given by Schur functions in the Heisenberg
generators. We were motivated to consider the Schur functions for the following reasons. The imaginary root vectors were constructed in \cite{CP}, where it was shown that they could be defined by a certain functional equation in terms of the Heisenberg generators.  After  a suitable renormalization, this equation is the
same as the equation that expresses the complete symmetric functions in terms
of the power sums. In Section 4 of this paper we show that the imaginary root 
vectors generate a polynomial algebra over $\bz[q,q^{-1}]$, are group-like with respect to the comultiplication and 
are quasi-orthonormal with respect to the  Drinfeld form on the algebra. It is well-known \cite{Ma} that the complete symmetric
functions are also group-like and  orthonormal with respect to the standard Hopf algebra structure and inner product on the ring of symmetric functions, and that the Schur functions form an orthonormal $\bz$-basis for the ring. Thus, we are able to identify the imaginary subalgebra with the ring of symmetric functions to construct the crystal basis of the imaginary part.
The appearance of Schur functions was anticipated in \cite{L6} in a series of conjectures on
level 0 representations of quantum affine algebras.
In  future work we shall  make explicit the connection between our work and these conjectures.

\vskip12pt\noindent{\bf Acknowledgements.} We thank N. Jing for helpful
conversations.  J.B. also thanks I. Frenkel and I. Grojnowski
for many helpful conversations.

\section{The algebras $\bu$ and $\bu^+$.}

In this section, we recall certain facts about $\ag$ and $\frak{g}$ and 
their associated quantum groups that will be needed later. 

Throughout this paper $\frak{g}$ will denote a simply-laced,
finite-dimensional complex simple Lie algebra and $(a_{ij})_{i,j\in I},
\ I=\{1,\dots, n\},$ will denote its Cartan matrix.   Let $(a_{ij})_{i,j\in\hat I}, \ 
\hat I=I \cup \{0\}$, be the extended Cartan matrix of $\frak{g}$ and let $\ag$ be the corresponding affine Lie algebra. Let $R$
(resp. $R^+$) denote a set of  roots (resp.  positive roots) of $\frak{g}$
and let $\alpha_i$ ($i\in I$) be a set of simple roots. Let $Q$ be the root lattice of $\frak{g}$, let $P$ be the
weight lattice, and let $\omega_i\in P$ ($i\in I$) be the fundamental weights of
$\frak{g}$.  For $\omega\in P$, $\eta\in Q$, define an integer
$|\omega|\cdot|\eta|$ by extending bilinearly the assignment
$|\omega_i|\cdot|\alpha_j| =\delta_{ij}$.  Notice that $|\alpha_i|\cdot|\alpha_j|
=a_{ij}$.  The corresponding objects for $\ag$ are defined similarly and
we denote them by $\hat{\cal{R}}$, $\hat\cal{Q}$ and so on.  Let $\theta$ be the
highest root of $\frak{g}$. Then, it is well-known that the element
$\alpha_0+\theta =\delta\in\hat{\cal{R}}$ and that $|\delta|\cdot|\alpha_i| =0$ for all
$i\in\hat{I}$. Further,  the set of roots $\hat\car$ of $\ag$ is
given by $\hat\car =\hat\car^+\cup -\hat\car^+$, where
\begin{equation*} \hat\car^+ = \{\alpha +k\delta\,\mid\, k\ge 0, 
\alpha\in R^+\}\cup\{k\delta\,\mid\,k>0\}\cup\{-\alpha+k\delta\,\mid\, k>0,
\alpha\in R^+\}.\end{equation*}  Set \begin{align*}\car_>& =\{k\delta+\alpha \,\mid\,k\ge 0, \alpha\in
  R^+\},\\ \car_<&=\{k\delta-\alpha \,\mid\, k> 0, \alpha\in R^+\},\\ 
  \car_0&=\{k\delta\,\mid\,k>0\}\times I,\\ \car&
  =\car_>\cup\car_0\cup\car_< .\end{align*} We call $\car$ the set of
  positive roots (with multiplicity) of $\ag$. Given an element
  $\eta\in\hat{\cal{Q}}$, let $re(\eta)\in Q$ be such that $\eta
  =k\delta+re(\eta)$ for some $k\in\bz$. We call $re(\eta)$ the real part of
  $\eta$.

  Let $W$ and $\hat{W}$ be the Weyl groups of $\frak{g}$ and $\ag$,
  respectively. It is well-known that they are Coxeter groups generated
  by simple reflections $s_i$ for $i\in I$ and $s_i$ for $i\in\hat{I}$,
  respectively. The Weyl group $W$ acts on the root lattice $Q$ by
  extending $s_i(\alpha_j) = \alpha_j-a_{ij}\alpha_i$.   Then,  $\hat{W}$ is isomorphic to the semi-direct product $W\tilde\times Q$, under the map $s_i\to (s_i, 0)$ for $i\in I$,
  $s_0\to(s_\theta,\theta)$, where $s_\theta(\alpha_j)=\alpha_j-(|\theta|\cdot|\alpha_j|)\theta$. The extended Weyl group $\tilde{W}$ is defined to be the
  semi-direct product $W{\tilde{\times}} P$. For any $w\in W$
  we write $(w,0)$ for the corresponding element in $\tilde{W}$, and for
  $\omega\in\P$ we write $t_\omega$ for the element $(1,\omega)$.
  The affine Weyl
  group $\hat{W}$ is a normal subgroup of $\tilde{W}$, and the quotient $\cal{T}=\tilde{W}/\hat{W}$ is a finite group 
  isomorphic to a subgroup of the group of diagram automorphisms
  of $\ag$, i.e. the bijections $\tau:\hat I\to\hat I$ such that
  $a_{\tau(i)\tau(j)}=a_{ij}$ for all $i,j\in\hat I$. Moreover, there is
  an isomorphism of groups $\tilde W\cong\cal{T}\tilde{\times}\hat W$,
  where the semi-direct product is defined using the action of $\cal{T}$
  in $\hat{W}$ given by $\tau.s_i=s_{\tau(i)}\tau$ (see \cite{Bo}). If
  $w\in\tilde W$, a reduced expression for $w$ is an expression $w=\tau
  s_{i_1}s_{i_2}\ldots s_{i_m}$ with $\tau\in\cal{T}$,
  $i_1,i_2,\ldots,i_m\in\hat I$ and $m$ minimal; we define the length $l(w)$ 
  of $w$ to be $m$.  The element $t_\omega$ acts on $\hat\cal{Q}$ by extending
  $t_{\omega}(\alpha_i) =\alpha_i-(|\omega|\cdot |\alpha_i|)\delta$.

  For $i\in I$, let $\tau_i\in \cal{T}$ be such that $w_i'=\tau_i^{-1}
  t_{\omega_i}\in\hat{W}$. Then, by \cite[Lemma 3.1]{B1}, we know that
  $l(w_i's_i)=l(w_i')-1$ and that there exists $j$ such that
  $l(s_jw_j')=l(w_j')-1$ and $\tau_i(s_j)=s_0\tau_i$.  Define elements
  $w_1,w_2,\dots, w_{2n}$ by
\begin{align*} w_i &=\tau_i\tau_{i-1}\ldots\tau_1 
w_i'(\tau_i\tau_{i-1}\ldots\tau_1)^{-1},\\
  w_{n+i}& = \tau_n\tau_{n-1}\ldots
  \tau_{i+1}(\tau_i\tau_{i-1}\ldots\tau_1)^2 w_i'(\tau_n\tau_{n-1}\ldots
  \tau_{i+1})^{-1}(\tau_i\tau_{i-1}\ldots\tau_1)^{-2},\end{align*}
if $1\le i\le n$. Let $2\rho$ denote the sum of all the roots
  in $R^+$. It is well-known that $2\rho=2\sum_{i=1}^n\omega_i$. The following lemma is easily established by using standard results on Coxeter groups (see \cite{Bo}, for instance).
\begin{lem}{\label{trho}} 
\begin{enumerate}
\item[(i)] The element $t_{2\rho}\in\hat{W}$, and has length
$N=\sum_i l(t_{2\omega_i})$.
\item[(ii)] $t_{2\rho}=w_1w_2\ldots w_{2n}$. 
\item[(iii)] There exists a reduced expression 
  $s_{i_1}s_{i_2}\ldots s_{i_N}$ for $t_{2\rho}$ such that the expressions
  $s_{i_1}\ldots s_{i_{l(t_{\omega_1})}}$, $s_{i_{l(t_{\omega_1})+1}}\ldots
  s_{i_{l(t_{\omega_1})+l(t_{\omega_2})}}, \dots $ etc., are reduced expressions for
  $w_1$, $w_2,\dots $ etc.
\item[(iv)] Define a doubly infinite sequence
\begin{equation*} \boh =(\ldots, i_{-1}, i_{0}, i_{1},\ldots )\end{equation*}
  by setting $i_k=i_{k(\!\!\!\mod N)}$ for $k\in\bz$. Then, for any integers $m<p$, the product
  $s_{i_m}s_{i_{m+1}}\ldots s_{i_p}$ is reduced.
\item[(v)] We have
\begin{equation*}
\begin{split} \car_> = \{ \alpha_{i_0}, s_{i_0} (\alpha_{i_{-1}}), 
s_{i_0} s_{i_{-1}} (\alpha_{i_{-2}}), \dots \}, \\
  \car_< = \{ \alpha_{i_1}, s_{i_1} (\alpha_{i_{2}}), 
s_{i_1} s_{i_{2}} (\alpha_{i_3}), \dots \}.\ \ \ \qedsymbol
\end{split}
\end{equation*}

\end{enumerate}
\end{lem} From now on we fix a reduced expression for
$t_{2\rho}$ as in (iii) above. For $k\in\bz$, write $k= k_0 +rN$, with $|k_0|<
N$ and $k$, $k_0$ and $r$ either zero or of the same sign, and set
\begin{align} 
  \beta_{k}&=s_{i_0}s_{i_{-1}}\ldots s_{i_{k+1}}(\alpha_{i_k}) =
  t_{2\rho}^{-r}s_{i_0}s_{i_{-1}}\ldots
  s_{i_{k_0+1}}(\alpha_{i_{k_0}})\quad {\text{if $k\le 0$},}\\ 
  \beta_k&=s_{i_1}s_{i_2}\ldots s_{i_{k-1}}(\alpha_{i_k}) =
  t_{2\rho}^{r} s_{i_1}s_{i_2}\dots s_{i_{k_0-1}}(\alpha_{i_{k_0}})
  \quad {\text{if $k>0$}}.
\end{align} Define a total order on
  $\car$ by setting
\begin{equation} \label{rootorder}
  \beta_0 < \beta_{-1} < \beta_{-2} \dots < \delta^{(1)} < \dots <
  \delta^{(n)} < 2 \delta^{(1)} < \dots < \beta_3 < \beta_2 < \beta_1,
\end{equation}
where $k\delta^{(i)}$ denotes $(k\delta,i)\in\car_0$.

Let $q$ be an indeterminate, let $\bq(q)$ be the field of rational
functions in $q$ with rational coefficients, and let $\bz[q,q^{-1}]$ be
the ring of Laurent polynomials with integer coefficients. For
$r,m\in\bn$, $m\ge r$, define
\begin{equation*} 
[m]=\frac{q^m -q^{-m}}{q -q^{-1}},\ \ \ \ 
  [m]! =[m][m-1]\ldots [2][1],\ \ \ \ 
\left[\begin{matrix} m\\ 
  r\end{matrix}\right] = \frac{[m]!}{[r]![m-r]!}.
\end{equation*}
  Then $\left[\begin{matrix} m\\r\end{matrix}\right]\in\bz[q,q^{-1}]$
  for all $m\ge r\ge 0$. 

\begin{prop}{\label{defnbu}} There is a Hopf algebra $\bu$ over $\bq(q)$ which is generated as an algebra by elements $E_{\alpha_i}$, $F_{\alpha_i}$, $K_i^{{}\pm 1}$ ($i\in\hat I$), with the following defining relations:
\begin{align*} 
  K_iK_i^{-1}=K_i^{-1}K_i&=1,\ \ \ \ K_iK_j=K_jK_i,\\ 
  K_iE_{\alpha_j} K_i^{-1}&=q^{ a_{ij}}E_{\alpha_j},\\ 
K_iF_{\alpha_j} K_i^{-1}&=q^{-a_{ij}}F_{\alpha_j},\\
  [E_{\alpha_i}, F_{\alpha_j}
]&=\delta_{ij}\frac{K_i-K_i^{-1}}{q-q^{-1}},\\ 
  \sum_{r=0}^{1-a_{ij}}(-1)^r\left[\begin{matrix} 1-a_{ij}\\ 
  r\end{matrix}\right]
&(E_{\alpha_i})^rE_{\alpha_j}(E_{\alpha_i})^{1-a_{ij}-r}=0\ 
  \ \ \ \ \text{if $i\ne j$},\\
\sum_{r=0}^{1-a_{ij}}(-1)^r\left[\begin{matrix} 1-a_{ij}\\ 
  r\end{matrix}\right]
&(F_{\alpha_i})^rF_{\alpha_j}(F_{\alpha_i})^{1-a_{ij}-r}=0\ 
  \ \ \ \ \text{if $i\ne j$}.
\end{align*}
The comultiplication of $\bu$ is given on generators by
$$\Delta(E_{\alpha_i})=E_{\alpha_i}\ot 1+K_i\ot E_{\alpha_i},\ \ 
\Delta(F_{\alpha_i})=F_{\alpha_i}\ot K_i^{-1} + 1\ot F_{\alpha_i},\ \ 
\Delta(K_i)=K_i\ot K_i,$$
for $i\in\hat I$.\hfill\qedsymbol
\end{prop}

Let $\bu^+$ (resp. $\bu^-$, $\bu^0$)  be the $\bq(q)$-subalgebras of $\bu$ generated by the $E_{\alpha_i}$  (resp. $F_{\alpha_i}$,  $K_i^{\pm 1}$) for $i\in \hat{I}$. The following result is well-known, see \cite{L2} for instance. 

\begin{lem}{\label {butriangle}} 
 $\bu\isom \bu^-\otimes\bu^0\otimes\bu^+$ as $\bq(q)$-vector spaces.\hfill\qedsymbol
\end{lem}

\begin{defn} \begin{enumerate}
\item[(i)] Let $\sigma:\bu^+\to\bu^+$ denote the $\bq(q)$-algebra anti-automorphism obtained by extending $\sigma(E_{\alpha_i}) =E_{\alpha_i}.$
\vskip 6pt
\item[(ii)] Let $\omega:\bu\to \bu$ be the $\bq(q)$-algebra automorphism defined by extending
\begin{equation*} \omega(E_{\alpha_i}) =F_{\alpha_i}, \ \ \omega(F_{\alpha_i})=E_{\alpha_i},\ \ \omega({K_i}) =K_i^{-1}.\end{equation*}
\end{enumerate}
\end{defn}

 It is convenient to use the following notation:
\begin{equation*}E_{\alpha_i}^{(r)}=\frac{E_{\alpha_i}^r}{[r]!}.\end{equation*} The elements $F_{\alpha_i}^{(r)}$ are defined similarly.

Corresponding to each element $w\in\tilde{W}$ one can define an
  automorphism $T_w:\bu\to\bu$ as follows.  Let $T_i$ ($i\in\hat I$) be
  the $\bq(q)$-algebra automorphisms of $\bu$ defined as follows (see \cite {L2}):
  \begin{align}\label{braid}
    T_i(E_{\alpha_i}^{(m)})&=(-1)^mq^{-m(m-1)}F_{\alpha_i}^{(m)}K_i^m,\ \ 
    T_i(F_{\alpha_i}^{(m)})=(-1)^mq^{m(m-1)}K_i^{-m}E_{\alpha_i}^{(m)},\\ 
    T_i(E_{\alpha_j}^{(m)})&=\sum_{r=0}^{-ma_{ij}}(-1)^{r}q^{-r}
    E_{\alpha_i}^{(-ma_{ij}-r)}E_{\alpha_j}^{(m)}E_{\alpha_i}^{(r)}\ 
    \ \text{if $i\ne j$},\\ 
    T_i(F_{\alpha_j}^{(m)})&=\sum_{r=0}^{-ma_{ij}}(-1)^{r}q^{r}
    F_{\alpha_i}^{(r)}
    F_{\alpha_j}^{(m)}F_{\alpha_i}^{(-ma_{ij}-r)}\ \ {\text{if
        $i\ne j$}}.
\end{align}
The finite group $\Cal T$ acts as $\bq(q)$-Hopf algebra automorphisms of $\bu$ by
$$\tau.E_{\alpha_i} = E_{\tau({\alpha_i})},\ \ \tau.F_{\alpha_i} = F_{\tau({\alpha_i})},\ \
\tau.K_i=K_{\tau(i)}, \ \ \ \text{for all $i\in\hat I$}.$$ If $w\in W$
has a reduced expression $w=\tau s_{i_1}\ldots s_{i_m}$, let $T_w$ be
the $\bq(q)$-algebra automorphism of $\bu$ given by $T_w=\tau
T_{i_1}\ldots T_{i_m}$. Then, $T_w$ depends only on $w$, and not on the
choice of its reduced expression. Further, it is shown in
\cite{B1} that $T_{\omega_i}(E_{\alpha_j}) =E_{\alpha_j},$ $T_{\omega_i}(F_{\alpha_j}) =F_{\alpha_j}$, $T_{\omega_i}(K_j)=K_j$ if
$i\ne j$, $T_{\omega_i}(K_i)=K_iC^{-1}$.

This leads us to  another realization of $\bu$, due to \cite{Dr}, \cite{B1},
\cite{Jing}.
\begin{thm}{\label{newr}} There is an isomorphism 
  of $\bq(q)$-Hopf algebras from $\bu$ to the algebra with generators
  $x_{i,r}^{{}\pm{}}$ ($i\in I$, $r\in\bz$), $K_i^{{}\pm 1}$ ($i\in I$),
  $h_{i,r}$ ($i\in I$, $r\in \bz\backslash\{0\}$) and $C^{{}\pm{1/2}}$,
  and the following defining relations:
\begin{align*}
  C^{{}\pm{1/2}}\ &\text{are central,}\\ K_iK_i^{-1} = K_i^{-1}K_i
  =1,\;\; &C^{1/2}C^{-1/2} =C^{-1/2}C^{1/2} =1,\\ K_iK_j =K_jK_i,\;\;
  &K_ih_{j,r} =h_{j,r}K_i,\\  K_ix_{j,r}^\pm K_i^{-1} &= q^{{}\pm
    a_{ij}}x_{j,r}^{{}\pm{}},\\ 
  [h_{i,r},h_{j,s}]&=\delta_{r,-s}\frac1{r}[ra_{ij}]\frac{C^r-C^{-r}}
  {q-q^{-1}},\\ [h_{i,r} , x_{j,s}^{{}\pm{}}] &=
  \pm\frac1r[ra_{ij}]C^{{}\mp {|r|/2}}x_{j,r+s}^{{}\pm{}},\\ 
  x_{i,r+1}^{{}\pm{}}x_{j,s}^{{}\pm{}} -q^{{}\pm
    a_{ij}}x_{j,s}^{{}\pm{}}x_{i,r+1}^{{}\pm{}} &=q^{{}\pm
    a_{ij}}x_{i,r}^{{}\pm{}}x_{j,s+1}^{{}\pm{}}
  -x_{j,s+1}^{{}\pm{}}x_{i,r}^{{}\pm{}},\\ [x_{i,r}^+ ,
  x_{j,s}^-]=\delta_{i,j} & \frac{ C^{(r-s)/2}\psi_{i,r+s}^+ -
    C^{-(r-s)/2} \psi_{i,r+s}^-}{q - q^{-1}},\\ 
\sum_{\pi\in\Sigma_m}\sum_{k=0}^m(-1)^k\left[\begin{matrix}m\\k\end{matrix}
\right]
  x_{i, r_{\pi(1)}}^{{}\pm{}}\ldots x_{i,r_{\pi(k)}}^{{}\pm{}} &
  x_{j,s}^{{}\pm{}} x_{i, r_{\pi(k+1)}}^{{}\pm{}}\ldots
  x_{i,r_{\pi(m)}}^{{}\pm{}} =0,\ \ \text{if $i\ne j$},
\end{align*}
for all sequences of integers $r_1,\ldots, r_m$, where $m =1-a_{ij}$, $\Sigma_m$ is the symmetric group on $m$ letters, and the $\psi_{i,r}^{{}\pm{}}$ are determined by equating powers of $u$ in the formal power series 
$$\sum_{r=0}^{\infty}\psi_{i,\pm r}^{{}\pm{}}u^{{}\pm r} = K_i^{{}\pm 1} exp\left(\pm(q-q^{-1})\sum_{s=1}^{\infty}h_{i,\pm s} u^{{}\pm s}\right).$$

The isomorphism is given by [B1]
$$x_{i,r}^+ =o(i)^rT_{\omega_i}^{-r}(E_{\alpha_i}), \ \ x_{i,r}^- =o(i)^rT_{\omega_i}^{r}(F_{\alpha_i}),$$
where $o:I\to\{\pm 1\}$ is a map such that $o(i)=-o(j)$ whenever $a_{ij}<0$ (it is clear that there are exactly two possible choices for $o$).
\hfill\qedsymbol\end{thm}

The following lemma is easily checked.
\begin{lem}{\label {defnOmega}} There exists a $\bq$-algebra anti-automorphism 
$\Omega$ of $\bu$ such that $\Omega(q)=q^{-1}$ and, for all $i\in I, r\in\bz$,
\begin{equation*}\Omega(x_{i,r}^\pm)=x_{i,r}^\pm,\ \ 
\Omega(\psi_{i,r}^\pm)=\psi_{i,r}^\pm,\ \ 
\Omega(h_{i,r}) =-h_{i,r},\ \  \Omega(C^{1/2})=C^{1/2}.\ \ \ \ \ \ \ \ \ \ \ \qedsymbol
\end{equation*}
\end{lem}

We now define a set of root vectors for each element of $\car_>\cup \car_<$.
\begin{align} \label{realroots}
E_{\beta_k}  = 
\begin{cases}    T_{i_0}^{-1} T_{i_{-1}}^{-1} \dots
  T_{i_{k+1}}^{-1} (E_{\alpha_{i_k}})
 & \text{if $k\le 0$}, \\
  T_{i_1} T_{i_{2}} \dots  T_{i_{k-1}} 
(E_{\alpha_{i_k}}) &\text{if $k>0$}.
\end{cases}
\end{align}
It follows from \cite[Proposition 40.1.3]{L2} that the elements $E_{\beta_k}\in\bu^+$. The elements $F_{\beta_k}\in \bu^-$ are defined similarly by replacing $E_{\alpha_i}$  by $F_{\alpha_i}$. The set 
\begin{equation*}
\{E_{\beta_k}, F_{\beta_k} \,\mid\,\beta_k\in\car_<\cup\car_>\}
\end{equation*}  
is called the set of real root vectors.  The next lemma is obvious from the definition of the root vectors
(recall that $l(T_{t_{2\rho}}) =N$).
\begin{lem}{\label{Trho}} $T_{t_{2\rho}}^{-1}(E_{\beta_k})\in \bu^+$ if and only if $\beta_k<\beta_N$. If $\beta_N\le \beta_k\le \beta_1$, then $T_{t_{2\rho}}^{-1}(E_{\beta_k})\in\bu^-\bu^0$.\hfill\qedsymbol
\end{lem}

 For $k>0$, $i\in I$, set 
\begin{equation*}
\tilde\psi_{k,i} =E_{k\delta -\alpha_i}E_{\alpha_i} -q^{-2}E_{\alpha_i}E_{k\delta-\alpha_i}
\end{equation*} 
and define elements $E_{k\delta,i}\in\bu^+$ by the functional equation
\begin{equation*} 
\exp\left((q-q^{-1})\sum_{k=1}^\infty \ed[k]{i} u^k\right) = 1 +
\sum_{k=1}^\infty  (q-q^{-1}) {\tilde{\psi}}_{k,i}u^k.
\end{equation*}

The next result relates the generators in Theorem \ref{newr}    to the root vectors defined above.  It is an easy consequence of the choice of the reduced expression for $t_{2\rho}$ (see \cite{B2}). \begin{lem}{\label{reldrin}} \begin{align*}  x_{i,k}^+  &
  = \begin{cases} o(i)^k E_{k\delta+\alpha_i} & \text{if $k\ge 0$},\\ 
  - o(i)^k  F_{- k \delta- \alpha_i } K_i^{-1} C^{k} &\text{if
    $k<0$},\end{cases} \\ x_{i,k}^- & =
  \begin{cases} - o(i)^k C^{-k} K_i E_{ k\delta-\alpha_i } &\text{if
      $k> 0$},\\ o(i)^k F_{ - k\delta+\alpha_i} & \text{if $k\le
      0$}. \end{cases} \end{align*}

Further, for $k>0$, $i\in I$, we have
\begin{align*} 
h_{i,k}&=
     o(i)^kC^{-k/2} E_{k\delta,i} \ \ \ \  \text{if $k> 0$,}\\
      \psi_{i,k}&=o(i)^k(q-q^{-1})C^{-k/2}K_i\tilde{\psi}_{k,i} \ \ \ \ \text{if
        $k>0$}.\ \ \ \ \ \ \ \qedsymbol
\end{align*} 

\end{lem}

The following is now an immediate consequence of Theorem \ref{newr} and Lemma \ref{reldrin}. 
\begin{prop} {\label{newrealization} }
For $i,j\in I$, $k,l\ge 0$, we have
\begin{align*} 
[E_{k\delta,i}, E_{l\delta,j}]& =0, \ \ \ \ {\text{if $k,l>0$}},\\
[E_{k\delta,i}, E_{l\delta\pm\alpha_j}] =\pm o(i)^ko(j)^k & \frac{1}{k}[ka_{ij}]E_{(k+l)\delta\pm\alpha_i},\ \ \ \  {\text{if $k>0, l\delta\pm\alpha_i\in\car$}},\\
E_{(k+1)\delta\pm\alpha_i}E_{l\delta\pm\alpha_j}-q^{\pm a_{ij}}E_{l\delta\pm\alpha_j}E_{(k+1)\delta\pm\alpha_i}&=\\ 
q^{\pm a_{ij}}E_{(l+1)\delta\pm\alpha_i}E_{k\delta\pm\alpha_j} & -E_{k\delta\pm\alpha_j}E_{(l+1)\delta\pm\alpha_i}, \ \ 
{\text{if $k\delta\pm\alpha_i, l\delta\pm\alpha_j\in\car$}},\\
\tilde\psi_{k+l,i} =E_{k\delta -\alpha_i}E_{l\delta+\alpha_i} & -q^{-2}E_{l\delta+\alpha_i}E_{k\delta-\alpha_i}, \ \ {\text{if $k,l>0$}}.\ \ \ \  \qedsymbol
 \end{align*}
\end{prop}

\begin{cor}{\label{tomega}} For $k>0$, $i,j\in I$,  we have
\begin{align*} 
T_{\omega_j}(E_{k\delta,i}) &= E_{k\delta,i},\\
T_{\omega_j}(E_{k\delta\pm\alpha_i})&=E_{k\delta\pm\alpha_i}\ \ 
{\text{if $i\ne j$,}}\\
T_{\omega_i}(E_{k\delta+\alpha_i})= E_{(k-1)\delta+\alpha_i},\ \ 
& T_{\omega_i}(E_{k\delta-\alpha_i})=E_{(k+1)\delta-\alpha_i},\ \ {\text{if $k>0$}.}\ \ \ \ \ \qedsymbol
\end{align*}
\end{cor}

Let $\bu^+(>)$ (resp. $\bu^+(<)$, $\bu^+(0)$) be the $\bq(q)$-subalgebra of $\bu^+$ generated by the $E_{\beta_k}$ for $k\le 0$ (resp. $E_{\beta_k}$ for $k>0$, $E_{k\delta, i}$ for $k>0$). We now define a basis for these algebras.

Following \cite[Section 3]{CP}, we introduce elements ${P_{k,i}}$ ($i\in I$, $k\in\bz$, $k\ge 0$) in $\bu^+(0)$ by $P_{0,i}=1$ and
\begin{equation*} P_{k,i} = -\frac{1}{[k]} \sum_{r=1}^k q^{k-r} \tilde
  {\psi}_{r,i}P_{k-r,i}.\end{equation*}
\begin{rem}
The $P_{k,i}$ defined here are actually $q^{-k}o(i)^kC^{k/2}$ times the $P_{i,k}$
in \cite{CP}.
\end{rem}
Setting $\tilde{P}_{k,i}=\Omega(P_{k,i})$, we find that
\begin{equation*} \tilde P_{k,i} = \frac{1}{[k]}
  \sum_{r=1}^k q^{r-k} \tilde{\psi}_{r,i}\tilde P_{k-r,i}.
\end{equation*}

It is proved in \cite[Section 3]{CP} that $P_{k,i}$ and $\tilde{P}_{k,i}$ satisfy  the
functional equations
\begin{equation} \label{integralimaginary}
  \sum_{k \ge 0}{ P_{k,i}} u^k = \exp \left( \sum_{k = 1}^\infty
  \frac{-E_{k\delta,i} u^k}{[k]}\right),
\end{equation} 
\begin{equation} \label{tildeP} \sum_{k \ge 0}{\tilde{P}_{k,i}} u^k = 
  \exp \left( \sum_{k = 1}^\infty \frac{E_{k\delta,i}
    u^k}{[k]}\right).
\end{equation} 
It is clear from the above that $\bu^+(0)$ is generated as a $\bq(q)$-algebra  by the $\tilde{P}_{k,i}$ for $k>0$, $i\in I$.

\begin{defn}
\begin{enumerate}
\item[(1)] Let $\bn^{\car}$ (resp. $\bn^{\car_>}$, $\bn^{\car_<}$,
$\bn^{\car_0}$) be the set of finitely supported maps from $\car$
(resp. $\car_>$, ${\car_<}$, ${\car_0}$) to $\bn.$ 
\item[(2)] For $\boc, \boc' \in \bn^{\car}$, we say that $\boc < \boc'$ if 
there exists $\beta\in\car$ such that $\boc(\beta)< \boc'(\beta)$ and 
$\boc(\beta') =\boc'(\beta')$ if $\beta'>\beta$.
\end{enumerate}\end{defn}

From now on, we shall identify $\alpha$ with the 
indicator function $1_\alpha \in \bn^{\car}$. 
\begin{defn} Given $\boc \in \bn^{\car}$,  define $E_{\boc}$ (resp. $E^\prime_{\boc}$) to be the
  monomial formed by multiplying the $(E_{\beta_k})^{(\boc(\beta_{k}))}$, for real
  roots $\beta_k$ ($k \in \bz$), and the $\tilde{P}_{k,i}$ ($k > 0$,
  $i\in I$), (resp. $\tilde{\psi}_{k,i}$, $k>0$, $i\in I$) in the order \eqref{rootorder}. \end{defn}

\begin{defn} Set \begin{align*}B &=\{E_\boc\,\mid\, \boc\in\bn^\car\},\\
B_> &=\{E_\boc\,\mid\, \boc\in\bn^{\car_>}\},\\
B_0&=\{E_\boc\,\mid\, \boc\in\bn^{\car_0}\},\\
B_0'&=\{E_\boc'\,\mid\, \boc\in\bn^{\car_0}\},\\
B_<&=\{E_\boc\,\mid\, \boc\in\bn^{\car_<}\}.\end{align*}.
 \end{defn}

\begin{prop}  \label{bu+triangle} 

\begin{enumerate}
\item[(i)] $\bu^+(0)\bu^+(<)$ and $\bu^+(>)\bu^+(0)$ are $\bq(q)$-subalgebras of
  $\bu^+$.
\vskip 6pt
\item[(ii)] As $\bq(q)$-modules, we have $\bu^+ \isom \bu^+(>) \otimes
  \bu^+(0) \otimes \bu^+(<).$
\vskip 6pt
\item[(iii)] The set $B_> $ (resp. $B_<$, $B_0'$) is a $\bq(q)$-basis
  of $\bu^+(>)$ (resp. $\bu^+(<)$, $\bu^+(0)$).
\vskip 6pt
\item[(iv)] The set $B_0$ is a  $\bq(q)$-basis of $\bu^+(0)$.
\vskip6pt
\item[(v)]  The set $B$ is a  $\bq(q)$-basis of $\bu^+$.
\end{enumerate}
\end{prop}
\begin{pf} Parts (i) through (iii) were proved in \cite {B1}, \cite{B2}. Part (iv) can be now  deduced easily by using the definition of the $\tilde{P}_{k,i}$. Part (v) follows from parts (ii) through (iv).\end{pf}

An element $x\in\bu^+$ is said to have homogeneity $\sum_i  d_i\alpha_i\in\cal{Q}^+$ if $x$ is a  $\bq(q)$-linear combination of products of the $E_{\alpha_i}$ ($i\in\hat{I}$) in which $E_{\alpha_i}$ occurs $d_i$ times for all $i$. In that case, $x$ is said to be homogeneous, we let $|x|=\sum_{i\in\hat{I}}d_i\alpha_i$ denote its homogeneity, and we write $ht({x}) =\sum_{i\in\hat{I}}d_i$. Any $x\in\bu^+$ is a finite sum of homogeneous elements.

We end this section by recalling some results from \cite{L2} that will be used crucially in later sections.

For $i\in I$, there exists a unique $\bq(q)$-linear map
$r_i:\bu^+\to\bu^+$ given by $r_i(1)= 0, \ \ 
r_i(E_{\alpha_j})=\delta_{i,j}$ for $j\in I$, and satisfying $r_i(xy)=q^{|y|\cdot|\alpha_i|}
  r_i(x)y+xr_i(y)$ for all homogeneous $x,y\in\bu^+$ ([L3, 1.2.13]).  Similarly, there
exists a unique $\bq(q)$-linear map $_ir:\bu^+\to\bu^+$ given by
$_ir(1)= 0, \ \ _ir(E_{\alpha_j})=\delta_{i,j}$, and satisfying $_ir(xy)=
{_ir}(x)y+q^{|x|\cdot |\alpha_i|} x{} {_ir}(y)$ for all homogeneous $x,y\in\bu^+$.

The following is proved in \cite[Lemma 38.1.2, Lemma 38.1.5, Proposition 38.1.6]{L2}.

\begin{prop} {\label{ui}}
\begin{enumerate} 
\item[(i)]
 $\{x\in\bu^+ \,\mid\, T_i(x)\in\bu^+\} = \{x\in\bu^+\,\mid\,{_ir}(x)=0\}$.
\vskip 6pt
\item[(ii)] $\{x\in\bu^+ \,\mid\, T_i^{-1}(x)\in\bu^+\} = \{x\in\bu^+\,\mid\,{r_i}(x)=0\}$.
\vskip 6pt
\item[(iii)] Defining $\bu^+[i]$ and ${^\sigma}\bu^+[i]$ to be the sets in parts (i) and (ii), respectively, we have
\begin{equation*} \bu^+ =\bigoplus_{r\ge 0} E_{\alpha_i}^{(r)}\bu^+[i], \  \ \bu^+ =\bigoplus_{r\ge 0} {^\sigma}\bu^+[i]E_{\alpha_i}^{(r)}.\ \ \ \ \ \ \ \ \ \qedsymbol
\end{equation*}
\end{enumerate}
\end{prop}

\section{ The integral form and the integral PBW basis} 

Let $\A = \bz[q,q^{-1}].$ We define the {\em divided powers
  subalgebra} $\ua$ to be the $\A$-subalgebra of
$\bu$ generated by the $K_i^{\pm 1}$,  $E_{\alpha_i}^{(r)}$, $F_{\alpha_i}^{(r)}$,
 for $i\in\hat{I}$ and $r\ge 0$. The $\A$-subalgebras $\ua^\pm$ of $\ua$ are defined  in the obvious way. The main result in this section is:

\begin{thm}{\label{Abasis}} The set $B$ is an $\A$-basis of $\ua^+$.\end{thm}
The rest of this section is devoted to proving Theorem 2. We proceed as follows. First we show that   $B$ is contained in $\ua^+$. This allows us to define $\A$-subalgebras
$\ua^+(>)$, $\ua^+(<)$ and $\ua^+(0)$. Then we prove that $B_>$, $B_< $ and $B_0$  are $\A$-bases of these subalgebras.
Finally we prove a triangular decomposition $\ua^+=\ua^+(>)\ua^+(0)\ua^+(<)$. By Proposition {\ref{bu+triangle}}, this means that, as $\A$-modules, $\ua^+\cong\ua^+(>)\ot\ua^+(0)\ot\ua^+(<)$. The theorem follows.

For $i\in\hat{I}$, $r\ge 1$, $m\in\bz$, define elements
\begin{equation*}\genfrac{[}{]}{0pt}{}{K_i,m}{r} = \prod_{s=1}^r \frac{K_i q^{m-s+1} - K_i^{-1}
    q^{-m+ s-1}}{q^s - q^{-s}}.
\end{equation*}
The following lemma is well-known (see \cite[Corollary 3.1.9]{L2}).

\begin{lem}  For $r,s\in\bz$, $r,s\ge 0$, we have
\begin{equation*} E_{\alpha_i}^{(r)}F_{\alpha_i}^{(s)} =\sum_{t=0}^{min(r,s)}F_{\alpha_i}^{(s-t)}\genfrac{[}{]}{0pt}{}{K_i,2t-r-s}{t}E_{\alpha_i}^{(r-t)}.
\ \ \ \ \ \qedsymbol
\end{equation*}
\end{lem}  

Let $\ua^0$ be the $\A$-subalgebra of $\ua$ generated by the $K_i^{\pm 1}$, $\genfrac{[}{]}{0pt}{}{K_i,m}{r}$ for all $i\in\hat{I}$, $r\ge 1$, $m\in\bz$. The following is well-known and can be deduced from the preceding lemma and Lemma \ref{butriangle}.
\begin{lem}{\label{uatriangle}} As $\A$-modules, we have
$\ua\cong\ua^-\ot\ua^0\ot\ua^+$. In particular,  we  have $\ua^+ =\bu^+\cap\ua$.\hfill\qedsymbol\end{lem}

The formulas given in Section 1 show that the restriction of $T_i$ to $\ua$ defines an $\A$-algebra automorphism of $\ua$. This proves, in light of the preceding lemma,  that $E_{\beta_k}^{(r)}\in\ua^+$ for all $k\in\bz$, $r\ge 0$. Let $\ua^+(>)$ (resp. $\ua^+(<)$) be the $\A$-subalgebras generated by the $E_{\beta_k}^{(r)}$ for $k\le 0$ (resp. $k>0$), $r\ge 0$. 

We next show that $\tilde{P}_{k,i}\in\ua^+$. For this we need several results from \cite{CP} which we now recall.

Let $\xi$ be an indeterminate and form the polynomial algebra
$\bq(q)[\xi]$ over $\bq(q).$ For $r\in\bz$, $r>0$, set $$\xi^{(r)} = \frac{\xi^r}{[r]!}.$$ For $i\in I$, define the following
elements of the algebra $\bu[[u]]$, where $u$ is another indeterminate:
\begin{equation*} \cx_i^+(u) = \sum_{k=0}^\infty \e{k\delta + \alpha_i}
  u^k, \ \  \cx_i^-(u) = \sum_{k=0}^\infty \e{(k+1)\delta - \alpha_i} u^k.
\end{equation*}
Let $\cd_i^\pm : \bq(q)[\xi] \to\bu[[u]]$ be the $\bq(q)$-algebra
homomorphisms that take $\xi$ to $\cx_i^\pm(u).$ Writing
\begin{equation*} \cd_i^\pm = \sum_{k = 0}^\infty D_{k,i}^\pm u^k,
\end{equation*}
the fact that $\cd_i^\pm$ are homomorphisms is equivalent
to $D_{k,i}^\pm(1)=\delta_{k,0}$ and
\begin{equation*}
D_{k,i}^\pm(fg) = \sum_{m=0}^kD_{m,i}^\pm(f) D_{k-m,i}^\pm(g)
\end{equation*}
for all $f, g \in \bq(q)[\xi].$  The $\bq(q)$-linear maps 
$D_{k,i}^\pm : \bq(q)[\xi] \rightarrow \bu$ are uniquely
determined by this relation together with
\begin{equation*} D_{k,i}^+(\xi) = \e{k \delta +\alpha_i}, \ D_{k,i}^-(\xi) =
  \e{(k+1) \delta - \alpha_i}, \ D_{k,i}^\pm(1) = \delta_{k,0}.
\end{equation*}
It is clear from the definition that, for all $r\ge 1$,
\begin{equation}{\label{homo}}|D_{k,i}^+(\xi^{(r)})|= k\delta +r\alpha_i, \ \  |D_{k,i}^-(\xi^{(r)})| = (k+1)\delta-r\alpha_i.\end{equation}

\begin{prop} \label{dformula} We have
  \begin{align} \tag{i} & D_{k,i}^+ (\xi^{(r)}) 
    = \sum_{t=1}^{r-1}(-1)^{t+1} q^{t(r-1)} E_{\alpha_i}^{(t)}
    D_{k,i}^+(\xi^{(r-t)}) \\ & \hskip 1.5in + q^{r(r-1)}
    T_{\omega_i}^{-1} D_{k-r,i}^+(\xi^{(r)}). \notag \\ \tag{ii}
    & D_{k,i}^-(\xi^{(r)}) = \sum_{t=1}^{r-1} (-1)^{t+1} q^{t(r-1)}
    D_{k,i}^-(\xi^{(r-t)}) E_{\delta - \alpha_i}^{(t)} \\ & \hskip 1.5in +
    (-1)^rq^{r(r-1)} T_{\omega_i} D_{k-r,i}^-(\xi^{(r)}). \notag
  \end{align}
The second term on the right-hand side of these equations is omitted if $k<r$.
\ \qedsymbol\end{prop}
\begin{pf} For this proof only, we write $\bu$ as $\bu_q(\hat{\frak g})$ to make explicit its dependence on the underlying finite-dimensional complex simple Lie algebra $\frak g$. Using Theorem 1, it is easy to see that there is a unique homomorphism of $\bq(q)$-algebras $\mu_i:\bu_q(\hat{sl}_2)\to\bu_q(\hat{\frak g})$ such that
\begin{equation*}
\mu_i(x_{1,k}^\pm)=o(i)^kx_{i,k}^\pm,\ \ \ \ \mu_i(K_1)=K_i,\ \ \ \mu_i(C^{1/2})=C^{1/2}
\end{equation*}
for all $k\in\bz$ (we take $o(1)=1$ for $sl_2$). By checking on the generators in Theorem 1, it is easy to see that 
\begin{equation*}
T_{\omega_i}\circ\mu_i=\mu_i\circ T_{\omega_1}.
\end{equation*}
In our present notation, Proposition 4.2 in \cite{CP} states that
\begin{equation*}
D_{k,1}^+ (\xi^{(r)}) = \sum_{t=1}^{r-1}(-1)^{t+1} q^{t(r-1)} E_{\alpha_1}^{(t)}
D_{k,1}^+(\xi^{(r-t)}) + q^{r(r-1)}T_{\omega_1}^{-1} D_{k-r,1}^+(\xi^{(r)}).
\end{equation*}
(The automorphism $T$ used in   Proposition 4.2 in \cite{CP} is $T_{\omega_1}^{-1}$; the definition of $T$ that precedes the statement of 
Proposition 4.2 in \cite{CP} is incorrect.) Applying $\mu_i$ to both sides of this equation gives part (i) of the proposition. Part (ii) is proved similarly.
\end{pf}

\begin{cor} {\label {dform}}

 Let $r,k\ge 1$. Then,
$$D_{k,i}^+(\xi^{(r)})=\sum\mu(s_0,s_1,\ldots)(E_{\alpha_i})^{(s_0)}(E_{\delta+\alpha_i})^{(s_1)}\ldots,$$
where the sum is over those non-negative integers $s_0,s_1,\ldots$ such that $\sum_ts_t=r$ and $\sum_t ts_t=k$, and the coefficients $\mu(s_0,s_1,\ldots)\in\A$. In particular, the coefficient of $E_{k\delta+\alpha_i}^{(r)}$ in $D_{kr}^+(\xi^{(r)})$ is $q^{kr(r-1)}$. An analogous statement holds for $D_{k,i}^-(\xi^{(r)})$ and the coefficient of  $E_{k\delta+\alpha_i}^{(r)}$ in $D_{(k-1)r}^-(\xi^{(r)})$ is $(-1)^{(k-1)r}q^{(k-1)r(r-1)}$.\end{cor}
\begin{pf} The first part is immediate from the proposition. The second part follows by induction on $r$, noting that the term involving $E_{k\delta+\alpha_i}^{(r)}$ can arise only from the second term on the right-hand side of the formula for $D_{kr}^+(\xi^{(r)})$ (and similarly for the other case).\end{pf}

  Next, if $x\in\bu$, let $L_x:\bu\to\bu$ (resp. $R_x:\bu\to\bu$) be left
  (resp. right) multiplication by $x$, and define
$$\bold{D}_{k,i}= L_{q^{k} P_{k,i}}D_{0,i}^+ + L_{q^{k-1}
  P_{k-1,i}}D_{1,i}^++\cdots +L_{P_{0,i}}D_{k,i}^+,$$
$$\text{(resp. } \bold{\tilde{D}}_{k,i} = R_{ q^{-k} \tilde{P}_{k,i}}
D^+_{0,i} + R_{q^{1-k} \tilde{P}_{k-1,i}} D^+_{1,i}+ \cdots +
R_{\tilde{P}_{0,i}} D^+_{k,i}).$$

\begin{prop}{\label{greatidentity}}
   Let $r,s\in\bn$, $i\in I$. Then,
\begin{align} & \tag{i} E_{\alpha_i}^{(r)}E_{\delta-\alpha_i}^{(s)}
  =\sum_{t=0}^{\text{min}(r,s)}\sum_{m+k=t} q^{2rs - tr - ts +
    t}D_{m,i}^-(\xi^{(s-t)})\bold{D}_{k,i}(\xi^{(r-t)}), \\ & \tag{ii} 
  E_{\delta - \alpha_i}^{(s)}E_{\alpha_i}^{(r)} =
  \sum_{t=0}^{\text{min}(r,s)}\sum_{m+k=t} q^{rt + st - 2sr +
    t}\bold{\tilde{D}}_{k,i}(\xi^{(r-t)}) D_{m,i}^-(\xi^{(s-t)}).
\end{align}
\end{prop}

\begin{pf} Part (i) is a restatement in the current notation of Lemma
  5.1 in \cite{CP}. And part (ii) follows from part (i) by applying $\Omega.$
\end{pf}

\begin{cor} $\tilde{P}_{k,i}\in\ua^+$ for all $k>0$,  $i\in I$.\end{cor}
\begin{pf} The proof proceeds  by induction on $k$. If $k=0$ the result is obvious. Assume that $\tilde{P}_{m,i}\in\ua^+$ for all $m<k$. Then, taking $r=s=k$ in the identity above, we see using Corollary \ref{dform} and the induction hypothesis that 
\begin{equation*}E_{\delta - \alpha_i}^{(k)}E_{\alpha_i}^{(k)} =\tilde{P}_{k,i} +{\text {terms in $\ua^+$}}.\end{equation*}
Since the left-hand side is in $\ua^+$ the result follows.
\end{pf}

Let $\ua^+(>)$, $\ua^+(<)$ and $\ua^+(0)$ be the $\A$-subalgebras of $\ua^+$ generated by the sets $\{E_{\beta_k}^{(r)}\,\mid\, k\le 0, r\ge 1\}$,  $\{E_{\beta_k}^{(r)}\,\mid\, k>0, r\ge 1\}$ and $\{\tilde{P}_{k,i}\,\mid\, k>0, i\in I\}$, respectively.
The following lemma is now obvious from Proposition \ref{bu+triangle}(iv) and the preceding corollary.
\begin{lem}{\label{b0}} The set $B_0$ is an $\A$-basis for $\ua^+(0)$.
\hfill\qedsymbol\end{lem}

The next lemma can be found in \cite{L2} and is an obvious consequence of Proposition \ref{ui}(iii) and the fact that the maps ${_ir}$ and $r_i$ preserve $\ua^+$.
\begin{lem}{\label{uai}} For all $i\in\hat{I}$, we have
\begin{equation*}\ua^+ =\bigoplus_r E_{\alpha_i}^{(r)}\ua^+[i], \  \ \ua^+ =\bigoplus_r{} ^\sigma_\cal{A}\bu^+[i]E_{\alpha_i}^{(r)},\end{equation*}
where $\ua^+[i] =\{x\in\ua^+\,\mid\, {_ir}(x)= 0\}$ and $ ^\sigma_\cal{A}\bu^+[i]=\{x\in\ua^+\,\mid\,r_i(x)=0\}$.\hfill\qedsymbol\end{lem}

\begin{prop} Let $s_{j_1}s_{j_2}\ldots s_{j_k}$ be an arbitrary reduced expression in $\hat{W}$. Define an $\A$-subalgebra 
\begin{equation*} \ua^+_{j_1,j_2,\dots, j_k} =\{x\in\ua^+\,\mid\, T_{j_k}T_{j_{k-1}}\dots T_{j_1}(x)\in \ua^-\ua^0\}.\end{equation*}
Then, the elements $E_{\alpha_{j_1}}^{(r_1)}(T_{j_1}^{-1}E_{\alpha_{j_2}}^{(r_2)})\ldots (T_{j_1}^{-1}T_{j_2}^{-1}\dots T_{j_{k-1}}^{-1}E_{\alpha_{j_k}}^{(r_k)})$, for $r_1,\ldots,r_k\ge 0$, form an $\A$-basis of $\ua^+_{j_1,j_2,\dots ,j_k}$.
\end{prop}
\begin{pf}
For $1\le l\le k$, we have
 \begin{align*} T_{j_k}T_{j_{k-1}}\dots T_{j_1}\left(T_{j_1}^{-1}T_{j_2}^{-1}\dots T_{j_{l}}^{-1}E_{\alpha_{j_{l+1}}}^{(r_{l+1})}\right)& = T_{j_k}T_{j_{k-1}}\dots T_{j_{l+1}}E_{\alpha_{j_{l+1}}}^{(r_{l+1})}\\
= (-1)^{r_{l+1}}q^{-{r_{l+1}}({r_{l+1}}-1)} & T_{j_k}T_{j_{k-1}}\dots T_{j_{l+2}}F_{\alpha_{j_{l+1}}}^{({r_{l+1}})}K_{j_{l+1}}^{r_{l+1}}.
\end{align*} 
By \cite[Proposition 40.2.1]{L2},  the expression $s_{j_k}s_{j_{k-1}}\ldots s_{j_{l+2}}$ is reduced and so by 
\cite[Lemma 40.1.2]{L2} (working with $F's$ rather than $E's$) we see that the right-hand side of the above equation is in $\ua^-\ua^0$ for all $r_1,\ldots,r_k\ge 0$. Since $\ua^-\ua^0$ is an $\A$-subalgebra of $\ua$, it follows that   $E_{\alpha_{j_1}}^{(r_1)}(T_{j_1}^{-1}E_{\alpha_{j_1}}^{(r_2)})\ldots (T_{j_1}^{-1}T_{j_2}^{-1}\dots T_{j_{k-1}}^{-1}E_{\alpha_{j_k}}^{(r_k)})\in \ua^+_{j_1,j_2,\dots, j_k}$. 
Further, it is proved in \cite[Proposition 40.2.1]{L2} that these elements are linearly independent. So, it remains to prove that they span $\ua^+_{j_1,j_2,\dots ,j_k}$. 

Let $x\in\ua^+_{j_1,j_2,\dots, j_k}$ be homogeneous, and write it as a sum \begin{equation*}x=\sum_{r_1}E_{\alpha_{j_1}}^{(r_1)}x_{r_1}\end{equation*}
according to the decomposition in Lemma \ref{uai}. Since $T_{j_1}x_{r_1}\in{}_{\A}\bu^+$ by Proposition 1.4, write
\begin{equation*}
T_{j_1}x_{r_1}=\sum_{r_2} E_{\alpha_{j_2}}^{(r_2)}x_{r_1,r_2},
\end{equation*} 
again according to the decomposition in Lemma \ref{uai}. Repeating this process, we define for $1\le l\le k$ a set of elements $x_{r_1,r_2,\dots ,r_{l}}\in\ua^+[j_{l}]$ satisfying
\begin{equation*} T_{j_{l}}x_{r_1,r_2,\dots ,r_{l}}=\sum_{r_{{l+1}}}E_{\alpha_{j_{l+1}}}^{(r_{l+1})}x_{r_1, r_2,\dots,r_{l+1}}, \end{equation*} 
for $1\le l< k$. This gives
\begin{align*} T_{j_k}T_{j_{k-1}}\ldots T_{j_1}x &= \sum_{r_1}(T_{j_k}T_{j_{k-1}}\ldots T_{j_1}E_{\alpha_{j_1}}^{(r_1)})(T_{j_k}T_{j_{k-1}}\ldots T_{j_1}x_{r_1})\\
&=\sum_{r_1, r_2}(T_{j_k}T_{j_{k-1}}\ldots T_{j_1}E_{\alpha_{j_1}}^{(r_1)})(T_{j_k}T_{j_{k-1}}\ldots T_{j_2}E_{\alpha_{j_2}}^{(r_2)})(T_{j_k}T_{j_{k-1}}\ldots T_{j_2}x_{r_1, r_2})\\
&\dots\dots\dots\dots\dots\dots\dots\dots\dots\dots\dots\dots\dots\dots\dots\dots\dots\dots\dots \\
&=\sum_{r_1,r_2,\dots, r_k} (T_{j_k}T_{j_{k-1}}\ldots T_{j_1}E_{\alpha_{j_1}}^{(r_1)})(T_{j_k}T_{j_{k-1}}\ldots T_{j_2}E_{\alpha_{j_2}}^{(r_2)})\dots \\
&\ \ \ \ \ \ \ \ \ \ \  \ldots(T_{j_k}T_{j_{k-1}}E_{\alpha_{j_{k-1}}}^{(r_{k-1})})(T_{j_k}E_{\alpha_{j_{k}}}^{(r_k)})x_{r_1,r_2,\dots, r_k}.\end{align*}
The left-hand side of the above equation is in $\ua^-\ua^0$, and since 
$T_{j_k}T_{j_{k-1}}\ldots T_{j_l}E_{\alpha_{j_l}}^{(r_l)}$ is in $\ua^-\ua^0$, it follows from Lemma \ref{uatriangle} that $x_{r_1,r_2,\dots, r_k}\in\A$. Applying 
$(T_{j_k}T_{j_{k-1}}\ldots T_{j_1})^{-1}$, we now  get the statement of the proposition.
\end{pf}

\begin{prop} $B_>$ is an $\A$-basis of $\ua^+(>)$.\end{prop}
\begin{pf} Let $x\in\ua^+(>)$. Then, using the definition of the root vectors, we see that there exists an integer $l\le 0$ and $i_0,i_{-1},\ldots,i_l\in\hat{I}$ such that $x\in\ua^+_{i_0,i_{-1},\ldots, i_{l}}$. The result is now immediate from the preceding proposition.\end{pf}

One can prove similarly (working with $r_i$ and replacing $T_i$ by $T_i^{-1}$): 
\begin{prop} $B_<$ is an $\A$-basis of $\ua^+(<)$.\hfill\qedsymbol\end{prop}
We omit the details.

Summarizing, we have proved that $B$ is an $\A$-basis of $\ua^+(>)\ua^+(0)\ua^+(<)$.

It remains then to prove the triangular decomposition. We shall need a number of subalgebras and subspaces of $\ua$. We collect them in the following definition.
\begin{defn}
\begin{enumerate}
\item[(i)] For $i\in I$, let $\ua_i^+(\gg)$ (resp. $\ua_i^+(\ll)$) be the $\A$-subalgebra  of
  $\ua^+(>)$ generated by the elements $E_{k\delta+\alpha_i}^{(r)}$
  (resp. $E_{(k+1)\delta-\alpha_i}^{(r)}$) for 
  $k\ge 0 ,r\ge 0$.
\vskip 6pt
\item[(ii)] Let $\ua^+(\gg)$ (resp. $\ua^+(\ll)$) be the $\A$-subalgebras of $\ua^+(>)$ 
\item[] (resp. $\ua^+(<)$) generated by $\ua_i^+(\gg)$ 
(resp. $\ua_i^+(\ll)$)
for all $i\in I$.
\item[(iii)] The $\A$-subalgebra of $\ua$ generated by $\ua^+(\ll)$ and the $F_{\alpha_i}^{(r)}$ for all $i\in I$, $r\ge 0$ will be denoted by ${{_\ca{\tilde \bu}}}^+(\ll)$. 
Let ${{_\ca{\tilde \bu}}}^+$ be the $\A$-subalgebra of $\ua$ generated by $\ua^+(\gg)$ and ${{_\ca{\tilde \bu}}}^+(\ll)$. Set  $ _\ca\tilde\bu(\Delta)=\ua^+(\gg)\ua^+(0){{_\ca{\tilde \bu}}}^+(\ll)$.

\end{enumerate}
\end{defn}

 The following is an obvious consequence of Corollary \ref{dform}.
\begin{lem}{\label{duagg}} For all $k,r\ge 0$, $i\in I$,  we have \begin{equation*} D_{k,i}^{+}(\xi^{(r)})\in \ua_i^+(\gg), \ \ D_{k,i}^{-}(\xi^{(r)})\in \ua_i^+(\ll).\ \ \ \ \ \ \ \ \ \qedsymbol
\end{equation*}
\end{lem}
 
\begin{prop} {\label{neww}} For all $r\ge 0$,  the element $E_{\alpha_0}^{(r)}$ is in the $\A$-subalgebra generated by the $(x_{i,k}^\pm)^{(s)}$ for $i\in I$, $k,s\ge 0$. Hence, $\ua^+$ is an $\A$-subalgebra of  ${{_\ca{\tilde \bu}}}^+$.\end{prop}
\begin{pf} Assume first that $\frak{g}$ is not of type $E_8$. Then, we can choose
\begin{enumerate}
\item[(i)] $i_0\in I$  such that $|\omega_{i_0}|\cdot|\theta| =1$;
\item[(ii)] an element $w$ in the subgroup of $W$ generated by 
$\{s_j\,\mid\,j\in I, j\ne i_0\}$ satisfying $w(\alpha_{i_0})=\theta$. \end{enumerate}
Hence, $t_{\omega_{i_0}}ws_{i_0}(\alpha_{i_0}) =\alpha_0$. Setting $w' = t_{\omega_{i_0}}ws_{i_0}$, it follows that $l(w') = l(t_{\omega_{i_0}}w)-1$, and also from \cite{L2} that $T_{w'}E_{\alpha_{i_0}}^{(r)} = E_{\alpha_0}^{(r)}$. Further, it is easy to check that $l(t_{\omega_{i_0}}w) =l(t_{\omega_{i_0}}) +l(w)$, and so  we get 
\begin{equation*} T_{w'}E_{\alpha_{i_0}}^{(r)}= T_{{\omega_{i_0}}w}T_{i}^{-1}E_{\alpha_{i_0}}^{(r)} = T_{{\omega_{i_0}}}T_w T_{i}^{-1}E_{\alpha_{i_0}}^{(r)} .\end{equation*}
Now, the element $T_w T_{i}^{-1}E_{\alpha_{i_0}}^{(r)}$ is in the 
$\A$-subalgebra generated by the $E_{\alpha_i}^{(l)}, F_{\alpha_i}^{(m)}$ for $i\in I$, $l,m\ge 0$ (since $ws_i\in W$). The result now follows by using Theorem \ref{newr} and the fact that that $T_{\omega_i}(E_{\alpha_j}) =E_{\alpha_j},$ 
 $T_{\omega_i}(F_{\alpha_j}) =F_{\alpha_j},$ if $i\ne j$.

The case of $E_8$ is somewhat more complicated, since we can only choose an element $i_0$ such that $|\omega_{i_0}|\cdot|\theta| =2$.  However, the argument above can be modified by using \cite[Lemma 2.7]{L0}. We omit the details. 

The second statement in the lemma now follows from Lemma \ref{reldrin} and the fact that $\tilde{\ua^+}$ contains the generators of $\ua^+$.
\end{pf}

\begin{prop} \label{triangle2} We have:
\begin{enumerate}
\item[(i)]  ${{_\ca{\tilde\bu}}}^+={{_\ca{\tilde\bu}}}(\Delta)$;
\item[(ii)] $_\ca\tilde{\bu}^+(\ll) \subset \ua^+(<)\ua^-\ua^0$;
\item[(iii)] $\ua^+\isom \ua^+(\gg)\ua^+(0)\ua^+(<).$
\end{enumerate}
\end{prop}
We prove this proposition in the remainder of this section.

We need the following commutation relations. \begin{prop}  Let $i,j\in I$, $k,r\ge 0$. We have
\begin{align*}
  \tag{i}  & {\tilde{P}_{k,j}}E_{r\delta+\alpha_i}=\sum_{s=0}^k
  o(i)^so(j)^s  \frac{[a_{ji}][a_{ji}+1]\cdots
[a_{ji}+s-1]}{[s]!}    E_{(r+s) \delta +
\alpha_i}{\tilde{P}_{k-s,j}}; \\
 \tag{ii} &   E_{(r+1)\delta - \alpha_i} {\tilde{P}_{k,j}} =\sum_{s=0}^k
  o(i)^so(j)^s  \frac{[a_{ji}][a_{ji}+1]\cdots
[a_{ji}+s-1]}{[s]!} {\tilde{P}_{k-s,j}} E_{(r+s+1) \delta -\alpha_i}.\\
\end{align*}
\end{prop}
\begin{pf} This is the same as for Lemmas 3.3 and 3.4 in \cite{CP}.
\end{pf}

\begin{cor} \label{commute}\begin{enumerate}
\item[(i)] If $a_{ij} =0$,\begin{equation*}  \tilde{P}_{k,j}E_{r\delta+\alpha_i}^{(s)} = E_{r\delta+\alpha_i}^{(s)}\tilde{P}_{k,j}.\end{equation*}
\item[(ii)] If $a_{ji}=-1$,
\begin{equation*}
 \tilde{P}_{k,j}E_{r\delta+\alpha_i}^{(s)}=\sum_{m=0}^k q^{m(s-m)}
E_{r\delta+\alpha_i}^{(s-m)}E_{(r+1)\delta + \alpha_i}^{(m)}\tilde{P}_{k-m,j}.
\end{equation*}
\item[(iii)] Fix $i\in I$, let $x\in\ua^+_i(\gg)$ be homogeneous, and assume that  $\tilde{P}_{k,i}x\in$ {\hbox{$\ua^+(>)\ua^+(0)$.}} Then, there exist elements  $x_s\in\ua_i^+(\gg)$ of homogeneity $|x|+s\delta$  such that 
\begin{equation*}\tilde{P}_{k,i}x =\sum_s x_s\tilde{P}_{k-s,i}.\end{equation*}
\end{enumerate}
Analogous results hold involving the $E_{(r+1)\delta-\alpha_i}$.
\end{cor}
\begin{pf} Parts (i)  and (ii) follow from Proposition 2.8 by a direct computation using the relations in  Proposition \ref{newrealization}. If $i=j$,  a repeated application of part (i) of Proposition 2.8 implies that we can write
\begin{equation*}\tilde{P}_{k,i}x =\sum_s y_s\tilde{P}_{k-s,i}.\end{equation*}
for some $y_s$ of homogeneity $|y_s| =|x|+s\delta$ in the $\bq(q)$-subalgebra of $\bu^+$ generated by $\{E_{k\delta+\alpha_i}\,\mid\,k \ge 0\}$.
On the other hand, since $\tilde{P}_{k,i}x\in\ua^+(>)\ua^+(0)$,  by Lemma \ref{b0} we can write  
\begin{equation*} \tilde{P}_{k,i}x =\sum_{\boc \in \bn^{\car_0}} x_\boc E_\boc\end{equation*}
for some elements
$x_\boc\in\ua^+(>)$. Since $B_0$ is a basis of $ \bu^+(0)$, we can now equate coefficients to get the result. \end{pf}

\begin{lem}{\label{step1}} Let $i,j\in I, k,r,s\ge
0$. 
 \begin{enumerate} \item[(i)] For $0\le t\le k$, there exist $x_t\in \ua^+_j(\gg)$ (resp. $x_t\in\ua_j^+(\ll)$) of homogeneity $|x_t| =(s+k-t)\delta+r\alpha_j$ (resp. $|x_t|=(s+k+1-t)\delta -r\alpha_j $) such that 

\begin{equation*}
\tilde{P}_{k,i}E_{s\delta+\alpha_j}^{(r)}=\sum_tx_t\tilde{P}_{t,i}\ \ \ \ \text{ (resp. $E_{(s+1)\delta-\alpha_j}^{(r)}\tilde{P}_{k,i} =\sum_t \tilde{P}_{t,i}x_t$).}
\end{equation*}
\item[(ii)] The elements $D_{k,j}^-(\xi^{(s)})E_{\alpha_i}^{(r)}$ can be written as $\A$-linear combinations of products $xyz$, where  $x\in\ua_i^+(\gg)$, $y\in \ua_i^+(0)$, $z\in\ua_i^+(\ll)$ are homogeneous and  
\begin{equation*} 
 |x| =r'\alpha_i+l\delta, \ \ |y|=m\delta, \ \ |z| =p\delta-s'\alpha_i,
\end{equation*}
 where $r'\le r$, $s'\le s$ and $l+m+p =k+1$.
\item[(iii)] There exist elements $x_t$ of homogeneity $|x_t| = (k-t)\delta-r\alpha_j$ in the $\A$-subalgebra generated by $\ua_j^+(\ll)$ and the $F_{\alpha_j}^{(s)}$ for $s\ge 0$, such that 
\begin{equation*}
F_{\alpha_j}^{(r)}\tilde{P}_{k,i}=\sum_t \tilde{P}_{t,i}x_t.
\end{equation*}
\end{enumerate}

 \end{lem}
\begin{pf} If $i\ne j$ part (i) was  proved in Corollary
\ref{commute}(i), (ii), and part (ii) is obvious from the defining relations. 
Let  $i=j$. From Corollary \ref{tomega} we see that $T_{t_{\omega_i}}\tilde{P}_{k,i} = \tilde{P}_{k,i}$ for all $k\ge 0,i\in I$.  Now, observe that (i) is equivalent (by applying $T_{t_{\omega_i}}^s$) to:
\vskip 6pt
\noindent (i$^\prime$)  There exist homogeneous elements $x_t\in \ua^+_i(\gg)$, with $|x| =(k-t)\delta+r\alpha_i$, such that 
$\tilde{P}_{k,i}E_{\alpha_j}^{(r)}=\sum_tx_t\tilde{P}_{t,i}.$
\vskip 6pt
 
We prove (i$^\prime$) and (ii) simultaneously by induction on $k$.

 If $k=0$, then (i$^\prime$) is obvious for all $r$ and (ii)  is clear from  Proposition \ref{greatidentity}, Lemma \ref{duagg} and \eqref{homo}. 

 Assume now that (i$^\prime$) and (ii) hold for all smaller values of $k$ and for all $r, s\ge 0$.  By a repeated application of (i$^\prime$), we can write, for  $m<k$ and any homogeneous $x\in\ua_i^+(\gg)$,
\begin{equation}{\label{px}}
\tilde{P}_{m,i}x=\sum_tx_t\tilde{P}_{m-t,i},\end{equation}
  where $x_t\in\ua^+_i(\gg)$ and $|x_t|=|x|+t\delta$. Next, note that
$E_{\delta-\alpha_i}^{(k)}E_{\alpha_i}^{(k)}E_{\alpha_i}^{(r)}$ belongs to $\ua^+_i(\gg)\ua_i^+(0)\ua^+_i(\ll)$  because it equals 
\begin{equation*}
\left[\begin{smallmatrix} {k+r} \\ k \end{smallmatrix} \right]
E_{\delta-\alpha_i}^{(k)}E_{\alpha_i}^{(r+k)}, 
\end{equation*}
 which is
in $\ua^+_i(\gg)\ua_i^+(0)\ua^+_i(\ll)$  by Proposition \ref{greatidentity}. On the
other hand, the same proposition implies that \begin{equation*}
E_{\delta-\alpha_i}^{(k)}E_{\alpha_i}^{(k)}E_{\alpha_i}^{(r)}=
\tilde{P}_{k,i}E_{\alpha_i}^{(r)} + wE_{\alpha_i}^{(r)}, \end{equation*} where
$w$ is a linear combination of terms of  type
$\bold{D}_{m,i}(\xi^{(l)})D^-_{t,i}(\xi^{(p)})$, with $m, t<k$. By equation \eqref{px} and the
induction hypothesis, we see that $wE_{\alpha_i}^{(r)}$ belongs to $\ua^+_i(\gg)\ua_i^+(0)\ua^+_i(\ll)$. Hence, $\tilde{P}_{k,i} E_{\alpha_i}^{(r)} \in \ua^+_i(\gg)\ua_i^+(0)\ua^+_i(\ll)$. Since $\tilde{P}_{k,i} E_{\alpha_i}^{(r)}\in\bu^+(>)\bu^+(0)$, we  conclude by Proposition \ref{bu+triangle} that $\tilde{P}_{k,i} E_{\alpha_i}^{(r)}\in$

\noindent$\ua_i^+(>)\ua^+(0)$.  But now applying Corollary \ref{commute}(iii),  we get (i$^\prime$) for $k$.

To prove (ii),  consider
\begin{equation} \label{twentysix} 
  E_{\delta-\alpha_i}^{(k+s)}E_{\alpha_i}^{(k)}E_{\alpha_i}^{(r)}=\left[
  \begin{smallmatrix} {k+r}
  \\ r
\end{smallmatrix}\right]E_{\delta-\alpha_i}^{(k+s)}E_{\alpha_i}^{(k+r)}.
\end{equation} By Proposition \ref{greatidentity}  and Lemma \ref{duagg}, the right-hand side
of \eqref{twentysix} belongs to 

\noindent $\ua^+_i(\gg)\ua_i^+(0)\ua^+_i(\ll)$ and the left-hand
side  equals \begin{equation*}
q^{-k(s-1)}D_{k,i}^-(\xi^{(s)})E_{\alpha_i}^{(r)}+wE_{\alpha_i}^{(r)},
\end{equation*} where $w$ is a linear combination of terms
$\bold{D}_{m,i}(\xi^{(l)})D^-_{t,i}(\xi^{(p)})$ with $m\le k$ and
$t<k$. The induction hypothesis and the fact that (i) holds for $k$ now
implies that $D_{k,j}^-(\xi^{(s)})E_{\alpha_i}^{(r)}\in \ua^+_i(\gg)\ua_i^+(0)\ua^+_i(\ll)\subset\ua^+(>)\ua^+(0)\ua^+(<)$.  To complete the proof of (ii), write \begin{equation*} D_{k,j}^-(\xi^{(s)})E_{\alpha_i}^{(r)}
=\sum_\boc a_\boc E_\boc,\end{equation*}
 where $a_\boc\in \A$. This is possible since we have already proved that $B$ is an $\A$-basis of $\ua^+(>)\ua^+(0)\ua^+(<)$.  Now, it is easy to see, using Proposition \ref{newrealization}, that 
$a_\boc \ne 0$ only if $\boc$ is supported on the set 
\begin{equation*}
\{l\delta+r\alpha_i\,\mid\,r'\le r\}\cup\{m\delta^{(i)}\,\mid\,m \le s\}\cup\{n\delta-s'\alpha_i\,\mid\,s'\le s\}.
\end{equation*}
The result follows. The case of $E_{s\delta-\alpha_j}^{(r)}\tilde{P}_{k,i}$ is similar.

Part (iii) can be deduced from part (i) by taking $s=0$ and applying $T_{t_{\omega_i}}^{-1}$ to part (i) and using Lemma \ref{reldrin}. \end{pf}

Since ${{_\ca{\tilde\bu}}}(\Delta)$ contains the generators of ${{_\ca{\tilde\bu}}}^+$, to prove part (i) of Proposition \ref{triangle2}, it suffices to show that ${{_\ca{\tilde\bu}}}(\Delta)$ is an $\A$-subalgebra of ${{_\ca{\tilde\bu}}}^+$. We prove that, if
$x\in\ua^+(\gg)$, $y\in\ua^+(0)$ and $z\in\tilde{\ua^+}(\ll)$ are homogeneous,
 then 
\begin{equation*} yx\in\ua^+(\gg)\ua^+(0), \ \ zy\in\ua^+(0){{_\ca{\tilde\bu}}}(\ll), \ \ zx\in{{_\ca{\tilde\bu}}}(\Delta).\end{equation*}
This is enough, in view of the following scheme (in which $x$'s (resp. $y$'s, $z$'s) denote elements of $\ua^+(\gg)$ (resp. $\ua^+(0)$, ${{_\ca{\tilde\bu}}}(\ll)$)):
\begin{align*}
(x_1y_1z_1)(x_2y_2z_2)&=\sum x_1(y_1x_3)y_3(z_3y_2)z_2,\ \ \ \ \text{if $z_1x_2=\sum x_3y_3z_3$}\\
&=\sum(x_1x_4)(y_4y_3y_5)(z_5z_2),\ \ \ \ \text{if $y_1x_3=\sum x_4y_4$, $z_3y_2=\sum y_5z_5$.}\end{align*}

To see that $yx\in\ua^+(\gg)\ua^+(0)$, we proceed by induction on $|y| =m\delta$. If $m=0$ there is nothing to prove. If $y=P_{m\delta,i}$ for some $ i$, the assertion is proved by using Lemma \ref{step1}(i) repeatedly. If  $y$ is  an $\A$-linear combination of terms of type $y_tP_{t\delta, j_t}$ for some $y_t\in\ua^+(0)$, $t\ge 0$, the result follows again by Lemma \ref{step1}(i) and the induction hypothesis. The case of $zy$ is similar.

To prove that $zx\in{{_\ca{\tilde\bu}}}(\Delta)$, we proceed by induction on the height of $re(-|z|)$. Consider first the case $ht(re(-|z|))=1$, so that for some $i\in I$ we have
\begin{equation*}
z=E_{l\delta-\alpha_i}\ \ \ \ \text{for some $l>0$,}\ \ \ \text{or}\ \ z=F_{\alpha_i}.
\end{equation*}
We can assume that $x=E_{m\delta+\alpha_j}^{(r)}$ for some $m,r\ge 0$, $j\in I$, for if $x$ is a product $E_{m_1\delta+\alpha_{j_1}}^{(r_1)}E_{m_2\delta+\alpha_{j_2}}^{(r_2)}\dots$ the proof can then be completed by an obvious induction along with the fact that we have proved that $\ua^+(\gg)\ua^+(0)$ is an $\A$-subalgebra of ${{_\ca{\tilde\bu}}}^+$.
Moreover, we can assume that $i=j$, since otherwise $z$ and $x$ commute. 

Consider first the case in which $z=E_{l\delta-\alpha_i}$ with $l>0$. Taking $s=1$, $k=l-1$ in  Lemma \ref{step1}(ii),   we see that $zx\in{{_\ca{\tilde\bu}}}(\Delta)$, and this proves the case $m=0$. To deal with the case $m>0$, we start by
using the fact that $B$ is an $\A$-basis of $\ua^+(>)\ua^+(0)\ua^+(<)$ to write $E_{l\delta-\alpha_i}E_{\alpha_i}^{(r)}$, $l>0$, as a sum
\begin{equation*} \sum_{\boc}a_\boc E_\boc\end{equation*}
where $a_\boc\in\A$. 
It is easy to see, by using Lemma \ref{newrealization} and an obvious induction on $r$, that
$a_\boc \ne 0$ only if the restriction of $\boc$ to $\bn^{\car_<}$ is either zero or supported on the root $l\delta-\alpha_i$. The result for arbitrary $m>0$  now follows by applying $T_{\omega_i}^{-m}$ to $E_{(l+m)\delta-\alpha_i}E_{\alpha_i}^{(r)}$ and making use of Corollary 1.1. 

If now $z=F_{\alpha_i}$, the case $m=0$ is contained in \cite[Corollary 3.1.9]{L2}. The case $m>0$ is deduced from this by applying 
$T_{\omega_i}^{-m}$ to $E_{m\delta-\alpha_i}E_{\alpha_i}^{(r)}$.

This completes the proof that $zx\in{{_\ca{\tilde\bu}}}(\Delta)$ when $ht(re(-|z|))=1$. In fact, more precisely we have proved that 
$zE_{m\delta+\alpha_i}^{(r)}$ can be writen as an $\A$-linear combination of products $x'z'$, where $x'\in\ua^+(\gg)\ua^+(0)$ and $z'\in{{_\ca{\tilde\bu}}}(\ll)$ are homogeneous and $|z'| =0$ or $|z'|=l'\delta-\alpha_i$ for some $l'\ge 0$.

Assume now that,  for all homogeneous elements $z$  such that $ht(re(-|z|))< s$ and all $x\in\ua^+(\gg)$,
$zx$ can be written as an $\A$-linear combination of terms $x'z'$ where $x'\in\ua^+(\gg)\ua^+(0)$ and $z'\in{{_\ca{\tilde\bu}}}(\Delta)$ has $ht(re(-|z'|))<ht(re(-| z|))$. We prove that the result holds for $ht(re(-|z|))=s$. If $z$ can be written as a product $z_1z_2$ such that $ht(re(-|z_t|))<ht(re(-| z|))$, for $t=1,2$, then we are done by induction. Otherwise, $z=E_{l\delta-\alpha_i}^{(p)}$ for some $i\in I, p\ge 0$. Arguing as in the case when $ht(re(-|z|))=1$, we see that it suffices to prove the result when $x=E_{\alpha_i}^{(r)}$ and $l>0$. 
Recall from Corollary \ref{dform} that  there exists an integer $M$ such that \begin{equation}
D^-_{(l-1)s,i}(\xi^{(s)})= q^ME_{l\delta-\alpha_i}^{(s)}+z'', \end{equation}
where $z''$ is a linear combination of products $\e{l_1
\delta - \alpha_i}^{(s_1)}\e{l_2 \delta -
  \alpha_i}^{(s_2)}\ldots$ with $l_1,l_2,\ldots$ and
$s_1,s_2,\ldots$ being positive integers such that $s_1,s_2,\ldots<s$
and $l_1s_1+l_2s_2+\cdots=ls$. Hence, 
\begin{equation*} \e{l \delta -
\alpha_i}^{(s)}E_{\alpha_i}^{(r)} = q^{-M}
\left(\pm D_{(l-1)s,i}^-(\xi^{(s)})E_{\alpha_j}^{(r)} -
zE_{\alpha_j}^{(r)}\right).  
\end{equation*} 
The first term has the correct form by Lemma \ref{step1}(ii), and the
second term  can be written (by a repeated
application of the induction hypothesis, since $s_i+r<s$) as an $\A$-linear combination of products $x'z'$, where $x'\in\ua^+(\gg)\ua^+(0)$ and $z'\in{{_\ca{\tilde\bu}}}(\ll)$ has $|z'|=l'\delta-s'\alpha_i$ for some $l'\ge 0, s'\le s$.  This proves part (i) of Proposition \ref{triangle2}.

To prove part (ii), observe (from the definition of the root vectors) that $T_{t_{2\rho}}$ maps $\tilde{\ua^+}(\ll)$ into $\ua^+(<)$.
Hence, for  $x\in\tilde{\ua^+}(\ll)$, we have
\begin{equation*}
T_{t_{2\rho}}(x)=\sum_{\boc\in\bn^{\car_<}} a_{\boc}E_\boc
\end{equation*}
for some $a_\boc\in{\cal A}$. Now, from Lemma \ref{Trho} we see that $T_{t_{2\rho}}^{-1}(E_\boc)\in \ua^+(<)\ua^-\ua^0$. Part (ii) follows.
Part (iii) is now obvious.

The proof of Proposition 2.7 is complete.\hfill\qedsymbol

\section{A computation of inner products on $\ua(0)$.}

The canonical basis is characterized by its behavior with respect to a
symmetric bilinear form introduced in \cite{Ka}, following Drinfeld.  Since we have introduced the imaginary root vectors, a
prerequisite to the construction of a crystal basis is an understanding
of the behavior of the form on $\ua^+(0).$ We begin with certain
preliminary definitions.

Define an
algebra structure on $\ua^+\ot\ua^+$ by
\begin{equation*} (x_1\ot x_2)(y_1\ot y_2) 
  =q^{|x_2|\cdot|y_1|}x_1y_1\ot x_2y_2,\end{equation*} where $x_t$, $y_t$
  ($t=1,2$) are homogeneous.  Let $r:\bu^+\to\bu^+ \ot \bu^+$ be the
  $\bq(q)$-algebra homomorphism defined by extending
$$r(E_{\alpha_i})
=E_{\alpha_i}\ot 1+1\ot E_{\alpha_i},\ \ \ \ (i\in\hat{I}).$$ 
The
algebra $\bu^+$ has a unique symmetric bilinear form $(\ , \ ):\bu^+ \times
\bu^+ \rightarrow \bq(q)$ \cite[1.2.5]{L2}  satisfying $(1,1)=1$ and
$$(E_{\alpha_i},E_{\alpha_j})=\delta_{i,j}(1-q^{-2})^{-1},\ \ 
(x,yy')=(r(x),y\ot y'), \ \ (xx',y)=(x\ot x',r(y)),$$ 
where the form on
$\bu^+\otimes\bu^+$ is defined by $(x_1\ot y_1, x_2\ot
y_2)=(x_1,x_2)(y_1,y_2).$ 
 The form satisfies
\begin{equation}{\label{rinner}} 
  (E_{\alpha_i}y, x) = (1-q^{-2})^{-1}(y, {_ir}(x))\ \ (yE_{\alpha_i},
  x) = (1-q^{-2})^{-1}(y, r_i(x)).\end{equation}

  Let $\ba= \bq(q)\cap \bq[[q^{-1}]]$. The main result of this section
  is:

\begin{prop}{\label{specialinner}} 
For $i,j\in I$, $k, k' >0$, we have 
\begin{equation*}
(\tilde{P}_{k,i},\tilde{P}_{k',j}) 
= \delta_{k,k'}\delta_{i,j} \ \ \ \mod\ (q^{-1}\ba).\end{equation*}
\end{prop}

For $i\in I$, let ${_\ca\Lambda_i}$ denote the $\ca$-subalgebra of
$\bu^+$ generated by the $\tilde{P}_{k,i}$ for $k\ge 0$. Since $\ua^+(0)$ is commutative, it follows from  Lemma \ref{b0} that ${_\ca\Lambda_i}$ is the polynomial algebra generated by the $\tilde{P}_{k,i}$
and that
$$\ua^+(0)\isom{_\ca\Lambda_1}\ot{_\ca\Lambda_2}\ot\cdots\ot {_\ca\Lambda_n}.$$ 
Any $x\in\ua^+(0)$ can thus be written  as a finite sum of products $x(1)\dots x(n)$, where $x(i)\in {_\ca\Lambda_i}$, and we denote this by
$$x=\sum x(1)\dots x(n).$$
  
Let $\Delta_i:  {_\ca\Lambda_i}\to {_\ca\Lambda_i}\ot{_\ca\Lambda_i}$ be the $\ca$-algebra homomorphism obtained by extending
\begin{equation*}\Delta_i(\tilde{P}_{k,i}) =\sum_{s=0}^k \tilde{P}_{s,i}\ot \tilde{P}_{k-s,i}.\end{equation*}
  Here, the algebra structure on ${_\ca\Lambda_i}\ot {_\ca\Lambda_i}$ is the
  usual one, namely
$$(x\ot y)(x'\ot y')= xx'\ot yy'\ \ \text{for}\ x,x',y,y'\in {_\ca\Lambda_i}.$$
\begin{cor}{\label{orthogonal}} Let $i,j\in I$. Then:
\begin{enumerate}
\item[(i)] if $x\in {_\ca\Lambda_i}$ and $y=y(1)\dots y(n)\in\ua^+(0)$ is such that $y(j)\ne 1$ for some $j\ne i$, then 
$$(x,y) = 0\ \ \ \mod\ (q^{-1}\ba);$$
\item[(ii)] for $x,y\in\ua^+(0)$, we have
  \begin{equation*}(x,y)=\sum(x(1),y(1))(x(2),y(2))\dots
    (x(n),y(n))\ \ \ \ \mod\ (q^{-1}\ba).\end{equation*}
\end{enumerate}
\end{cor}

The proof of Proposition 3.1 and Corollary 3.1 occupies the rest of this section. We start with some preliminary results.

\begin{prop}{\label{ri}}
\begin{enumerate}
\item[(i)] For $m>0$ and $i\in I$, we have
\begin{align*} & r_i(E_{\delta-\alpha_i}) = 0, \  \quad
   r_i(E_{\delta,i})=(1-q^{-4})E_{\delta-\alpha_i},\\ &
  r_i(E_{2m\delta-\alpha_i})= q(1-q^{-2})(1-q^{-4}) \\ & \hskip 1in \times
  \bigl(\sum_{s=0}^{m-2}q^{2s}
  E_{(2m-s-1)\delta-\alpha_i}E_{(s+1)\delta-\alpha_i}+
  q^{2m-3}E_{m\delta-\alpha_i}^{(2)}\bigr),\\ &
  r_i(E_{(2m+1)\delta-\alpha_i}) =
  q(1-q^{-2})(1-q^{-4})\sum_{s=0}^{m-1}q^{2s}
  E_{(2m-s)\delta-\alpha_i}E_{(s+1)\delta-\alpha_i}.\end{align*}
\item [(ii)] For $k>0$, $i,j\in I$ and $a_{ij}=-1$, we have
\begin{align*} & {_ir}(E_{k\delta-\alpha_j})=0, \ \quad
{_ir}(\tilde{\psi}_{k,j})=0,\\
& {_ir}(E_{k\delta,j}) =0,\ \quad  r_i(\frac{k}{[k]}E_{k\delta,j})
=-q^{-1}(1-q^{-2})E_{k\delta-\alpha_i}.\end{align*}
\item[(iii)] For $k>0$, $i,j\in I$ and $a_{ij}=0$, we have
\begin{equation*} 
r_i(E_{k\delta-\alpha_j}) =0,\ \ r_i(E_{k\delta, j}) =0.\end{equation*}
\end{enumerate}
\end{prop}
\begin{pf} The first equality in (i) follows from 
  \cite[Lemma 3.4]{B1}. The second is now easily deduced from the
  definition of $r_i$ and the relation
\begin{equation*}
E_{\delta-\alpha_i}E_{\alpha_i}-q^{-2}
E_{\alpha_i}E_{\delta-\alpha_i}=E_{\delta,i}.\end{equation*}
  The third and fourth equalities are proved by induction using the
  definition of $r_i$, the relation
\begin{equation*}[E_{\delta,i},E_{k\delta-\alpha_i}]=
-[2]E_{(k+1)\delta-\alpha_i}, \end{equation*}
and  the following consequence of the relation in 
Proposition \ref{newrealization} between the $E_{k\delta-\alpha_i}$ for  $k>0$:
\begin{align*}
  &  E_{\delta-\alpha_i}E_{2m\delta-\alpha_i}
  -q^2E_{2m\delta-\alpha_i}E_{\delta-\alpha_i}
=(q^4-1)\sum_{s=1}^{m-1}q^{2s-2}
E_{(2m-s)\delta-\alpha_i}E_{(s+1)\delta-\alpha_i},\\ 
&  E_{\delta-\alpha_i}E_{(2m+1)\delta-\alpha_i}
  -q^2E_{(2m+1)\delta-\alpha_i}E_{\delta-\alpha_i} \\ & \hskip .75in =
  (q^{4}-1)\bigl(\sum_{s=0}^{m-1}q^{2s-2}
E_{(2m+1-s)\delta-\alpha_i}E_{(s+1)\delta-\alpha_i}
  + q^{2m-3}E_{(m+1)\delta-\alpha_i}^{(2)}\bigr).\end{align*} We omit
  the details.

  The proof of the first two equalities in (ii) is similar, but instead
  using that ${_ir}(E_{\delta-\alpha_j})=0$ if $i\ne j$ (\cite[Lemma 3.5]{B1}).
  To prove the third equality we proceed as follows.  Assume by
  induction that ${_ir}(E_{s\delta ,j})=0$ for all $s<k$.  It is easy to
  deduce from Lemma 1.5 and the functional equation in Theorem \ref{newr} relating the ${\psi}_{m,j}$ to the $h_{m,j}$
  that $E_{k\delta,j}$ is in the $\A$-subalgebra generated by
  $\tilde{\psi}_{k,j}$ and the $E_{s\delta,j}$ for $s<k$. The result is now
  clear.  To prove the fourth equality, recall that for $x\in\bu^+$ we
  have, by \cite[Proposition 3.6]{L2},
\begin{equation*} 
  [x, F_{\alpha_i}]
  =\frac{r_i(x)K_i-K_i^{-1}{_ir}(x)}{q-q^{-1}}.\end{equation*} 
The result now follows by taking
  $x=\frac{k}{[k]}E_{k\delta,j}$ and using the defining relation
\begin{equation*}[\frac{k}{[k]}E_{k\delta,j}, F_{\alpha_i}]= 
  -K_iE_{k\delta-\alpha_i}\end{equation*}
 in
  Theorem \ref{newr} and Proposition \ref{newrealization}. The proof of (iii) is similar to that of
  (i), but instead using that $r_i(E_{k\delta-\alpha_j})
  =0$ if $a_{ij}=0$ (\cite[Lemma 3.5]{B1}).
\end{pf}

The following result is proved in \cite[Proposition 40.2.4]{L2}.
\begin{prop}\label{luinner} For $\boc,\boc'\in N^{\car}$, we have
\begin{equation*} 
  (E_\boc,E_{\boc'}) = (E_{\boc_0},E_{\boc'_0}) \prod_{s\in\bz}
  (E_{\alpha_{i_s}}^{(c_s)}E_{\alpha_{i_s}}^{(c_s)}),
\end{equation*}
  where the $i_s$ are as in Lemma $1.1 (iii)$ and $\boc_0$ (resp. $\boc_0'$) denotes the restriction of $\boc$ (resp. $\boc'$) to $\car_0$.
\hfill\qedsymbol
 \end{prop}

\begin{lem}{\label{psipsi}} Let $k>0$, $i,j\in I$. 
\begin{enumerate}
\item[(i)] We have
 $$(\tilde{\psi}_{k,i}, \tilde{\psi}_{k,i}) = \frac{q^{2k-2}(1-q^{-4})}{(1-q^{-2})^2}.$$
\item[(ii)] If $a_{ij}= -1$, then
  \begin{equation*}(\tilde{\psi}_{k,i},\frac{k}{[k]}E_{k\delta,j})
    =\frac{q^{-1}}{(1-q^{-2})}.\end{equation*}
\item[(iii)] If $a_{ij}=0$, then $(\tilde{\psi}_{k,i}, \tilde{\psi}_{k,j})=0$. 
\end{enumerate}\end{lem}
\begin{pf} Using Proposition 1.1, we see  that
\begin{align*}  (\tilde{\psi}_{k,i}, \tilde{\psi}_{k,i})&=(E_{k\delta-\alpha_i}E_{\alpha_i}-q^{-2}E_{\alpha_i}E_{k\delta-\alpha_i}, E_{k\delta-\alpha_i}E_{\alpha_i}-q^{-2}E_{\alpha_i}E_{k\delta-\alpha_i} )\\
  &=(1-q^{-2})^{-1}(E_{k\delta-\alpha_i},
  r_i(E_{k\delta-\alpha_i}E_{\alpha_i}))\\&-
  2q^{-2}(1-q^{-2})^{-1}(E_{k\delta-\alpha_i},
  r_i(E_{\alpha_i}E_{k\delta-\alpha_i}))\\ &+ q^{-4}(1-q^{-2})^{-1}(E_{k\delta-\alpha_i}, E_{k\delta-\alpha_i})\\ 
  &=\frac{1-q^{-4}}{(1-q^{-2})^2} +
  \frac{q^2}{(1-q^{-2})^2}(r_i(E_{k\delta-\alpha_i}),
  r_i(E_{k\delta-\alpha_i})).
\end{align*} 
The second and third equalities
  in this computation follow from \eqref{rinner} and Proposition \ref{luinner},
keeping in mind that  $E_{\alpha_i}E_{m\delta-\alpha_i}\in B$ for all $m>0$, and using the explicit formula for $r_i(E_{k\delta-\alpha_i})$ given in Proposition 3.2. The computation can now be completed by using Proposition \ref{ri},
  Proposition \ref{luinner}, and the following identity in \cite[Lemma
  1.4.4]{L2}
\begin{equation*} 
  (E_{\alpha_i}^{(p)}, E_{\alpha_i}^{(p)}) =
  \prod_{s=1}^p(1-q^{-2s})^{-1}.\end{equation*}

  To prove (ii), notice that
\begin{align*}(\tilde{\psi}_{k,i}, \frac{k}{[k]}E_{k\delta, j})&=
  (E_{k\delta-\alpha_i}E_{\alpha_i}-q^{-2}E_{\alpha_i}E_{k\delta-\alpha_i},
  \frac{k}{[k]}E_{k\delta, j})\\ &=(E_{k\delta-\alpha_i}E_{\alpha_i},
  \frac{k}{[k]}E_{k\delta, j}), \ \ \text{since} \ 
  {_ir}(E_{k\delta, j}) =0\\ &= (1-q^{-2})^{-1}(E_{k\delta-\alpha_i},
  q^{-1}(1-q^{-2})E_{k\delta-\alpha_i})\\ 
  &=\frac{q^{-1}}{(1-q^{-2})},\end{align*} where the penultimate
  equality follows from Proposition \ref{ri}(ii).  The proof of (iii) is
  similar, using Proposition \ref{ri}(iii).
\end{pf}

 We need some additional results
about the behaviour of  $\Delta$ and $r$ on the Heisenberg generators.

Consider the polynomial ring $\bq[x_k, \ k > 0]$. This has a
natural Hopf algebra structure with comultiplication obtained by extending 
$ x_k\to x_k\ot 1+1\ot x_k$. 
Defining elements $\lambda_k$ for $k\ge 0$ by $\lambda_0=1$ and 
$\lambda_k =\frac1k\sum_{s=1}^{k}x_s\lambda_{k-s}$ for $k>0$, and setting $\Lambda=\sum_{k=0}^\infty \lambda_k u^k$, it is a result of
\cite[appendix]{Ga} that the comultiplication takes $\Lambda$ to $\Lambda\ot\Lambda$.   

For any algebra $A$, let $A_+$ denote the augmentation algebra.
\begin{prop}{\label{comult}} Let $k>0$, $i\in I$. We have:
  \begin{align} \tag{i} & r(E_{k\delta,i}) = E_{k\delta ,i}\ot 1 +1 \ot
  E_{k\delta, i}+ {\text{terms in }} \bu^+(<)\bu^+(0)\otimes
  \bu^+(0)\bu^+(>)_+; \\
 \tag{ii} & r(\tilde{\psi}_{k,i}) =
\tilde{\psi}_{k,i}\ot 1 +1\ot\tilde{\psi}_{k,i} \\& +
(q-q^{-1})\sum_{s=1}^{k-1}\tilde{\psi}_{s,i}\ot\tilde{\psi}_{k-s,i}
+{\text{terms in }} \bu^+(<)_+\bu^+(0)\otimes \bu^+(0)\bu^+(>)_+; \notag
\\ \tag{iii} & r(\tilde{P}_{k,i})
=\sum_{s=0}^{k} \tilde{P}_{s,i}\ot \tilde{P}_{k-s,i} + {\text{terms in
    }} \bu^+(<)_+\bu^+(0)\otimes \bu^+(0)\bu^+(>)_+.
\end{align}
\end{prop}
\begin{pf}  The following formula for
  the coproduct of the $E_{k \delta, i}$ is proved in \cite[Proposition 7.1]{Da} and can be derived from the
  explicit coproduct formulas for the loop-like generators (\cite{B1},
  Proposition  5.3): 
$$ \Delta\bigl(\e{k \delta,
    i}\bigr) = \e{k\delta,i} \ot 1 + K_{k\delta} \otimes \e{k\delta,i} +
  \text{terms in } \bu^0\bu^+(<)_+\bu^+(0) \otimes \bu^+(>)_+
.$$ 
Similar formuals can also be found in \cite{JKK}. Part (i) now follows from \cite[3.1.5]{L2}. Part (ii)
  can be deduced from (i) (again a proof can be found in \cite{JKK}).

To prove (iii), recall from Section 1 that 
\begin{equation*} 
  \tilde{P}_{k,i} =\frac{1}{k}\sum_{s=1}^{k}\frac{s}{[s]}
  E_{s\delta.i}\tilde{P}_{k-s,i}.\end{equation*} 
 The proof of (iii) is now easily deduced from Lemma \ref{b0} and 
  \eqref{integralimaginary} by replacing $x_k$ by $\frac{k}{[k]}E_{k\delta,i}$ and using
  the fact that $\bu^+(>)_+\bu^+(0)\subset\bu^+(0)\bu^+(>)_+$.
  \end{pf}

\begin{cor}
\begin{enumerate} \item[(i)] If $x\in {_\ca\Lambda_i}$, then \begin{align*}
  r(x)& = \Delta_i(x)+ {\text{terms in }} \bu^+(<)_+\bu^+(0)\otimes
  \bu^+(0)\bu^+(>)_+ \\ &= x\ot 1+1 \ot x + \text{terms in }
  ({_\ca\Lambda_i})_+\ot ( {_\ca\Lambda_i})_+\\ &\ \ \ \ \ \ +{\text{terms in }}
  \bu^+(<)_+\bu^+(0)\otimes \bu^+(0)\bu^+(>)_+. \end{align*}
\item[(ii)] If $x\in {_\ca\Lambda_i}$, $y,z\in\ua^+(0)$, then
\begin{equation*}(x,yz)= (\Delta_i(x), y\ot z).\end{equation*}
\end{enumerate}
\end{cor}
\begin{pf} It suffices to prove (i) when $x$ is a product of the
$\tilde{P}_{k,i}$ for $k>0$, $i\in I$. If
$x$ has length one this is Proposition 3.4(iii). The proof can be
completed easily by an induction on the length of $x$, keeping in mind
that for $x,y\in\ua^+(0)$ we have $|x| \cdot |y|=0$. For part (ii), we have, by
Proposition \ref{luinner}, that 
$$(\bu^+(0)\bu^+(>)_+\,,\, \bu^+(0)) =0.$$
The result now follows from part (i), since $(x,yz)=(r(x),y\ot z)$ for
all $x,y,z\in\ua^+$.
\end{pf}
\begin{lem}{\label{psip}} For $k>0$, we have
\begin{enumerate}
\item[(i)] \begin{equation*} (\tilde{\psi}_{k,i}, \tilde{P}_{k,i})
  =\frac{q^{k-1}(1-q^{-2k-2})}{(1-q^{-2})^2};\end{equation*}
\item[(ii)] if $i\ne j$,
\begin{equation*}
(\tilde{\psi}_{k,i}, \tilde{P}_{k,j}) =0  \ \  \mod (   q^{-1}\ba).\end{equation*}
\end{enumerate}
\end{lem}
\begin{pf} The proof of (i) is by induction on $k$. The case $k=1$ is contained
  in Lemma \ref{psipsi}. Assume the result for all smaller values of $k$. Using the definition of $\tilde{P}_{k,i}$ we get
\begin{align}(\tilde{\psi}_{k,i}, \tilde{P}_{k,i})&=
  (r(\tilde{\psi}_{k,i}), \frac{q-q^{-1}}{[k]}\sum_{s=1}^{k-1}
  q^{s-k}\tilde{\psi}_{s,i}\ot \tilde{P}_{k-s,i}) +
  (\tilde{\psi}_{k,i},\tilde{\psi}_{k,i})\\ 
  &=\frac{q-q^{-1}}{[k]}\sum_{s=1}^{k-1}(q^{s-k}\tilde{\psi}_{s,i},
  \tilde{\psi}_{s,i})(\tilde{\psi}_{k-s,i}, \tilde{P}_{k-s,i})+
  (\tilde{\psi}_{k,i},\tilde{\psi}_{k,i}) ,\end{align} where we use
  Proposition \ref{comult} to get the last equality.  A direct
  computation using Lemma \ref{psipsi} and the induction hypothesis gives (i).

The proof of (ii) is identical except that we use the relation
\begin{equation*}
  \tilde{P}_{k,i}= \frac{1}{k}\sum_{s=1}^{k}
\frac{s}{[s]}E_{s\delta.i}\tilde{P}_{k-s,i}.\end{equation*}
  We omit the details.
\end{pf}

We can now complete the proof of Proposition \ref{specialinner} as
follows. If $k\ne k'$ the result is obvious since the elements have
different homogeneity. If $k=k'=1$ the result follows from Lemma
\ref{psipsi}. Assume the result for all smaller values of $k=k'$.  By using the definition of $\tilde{P}_{k,i}$,  the properties of the inner product and Proposition
\ref{comult}, we get
\begin{align*}
  (\tilde{P}_{k,i},\tilde{P}_{k,j})&=\frac{1}{[k]}\left(\sum_{s=1}^kq^{s-k}\tilde{\psi}_{s,i}\ot\tilde{P}_{k-s,i}
  \ ,\ \sum_{s=0}^k\tilde{P}_{s,j}\ot\tilde{P}_{k-s,j}\right)\\ 
  &=\frac{1}{q^{k-1}(1+q^{-2}+\dots
    +q^{-2k+2})}\sum_{s=1}^kq^{s-k}(\tilde{\psi}_{s,i},\tilde{P}_{s,j})(\tilde{P}_{k-s,i},\tilde{P}_{k-s,j})\\ 
  &=\frac{1}{(1+q^{-2}+\dots
    +q^{-2k+2})}\sum_{r=1}^kq^{s-2k+1}(\tilde{\psi}_{s,i},\tilde{P}_{s,j})(\tilde{P}_{k-s,i},\tilde{P}_{k-s,j}).\end{align*}
  Since $s-2k+1<0$, it is obvious by Lemma \ref{psip}(ii) that the right-hand side is in $q^{-1}\ba$ if $i\ne j$. If $i=j$, then by Lemma
  \ref{psip}(i) we see that the right-hand side is
\begin{equation*}
  \frac{1-q^{-2k-2}}{(1-q^{-2})^2(1+q^{-2}+\dots
    +q^{-2k+2})}\sum_{s=1}^kq^{s-k}(\tilde{P}_{k-s,i},\tilde{P}_{k-s,i})
  = 1 \ \ \  \mod (   q^{-1}\ba), \end{equation*} and the proof is complete.
  \hfill\qedsymbol

We turn now to the proof of Corollary \ref{orthogonal}. To prove (i), we
set $y=y(j)y'$ and write $\Delta_i(x)=\sum_sz_s\ot z_s'$, $z_s,z_s'\in{}_\ca\Lambda_i$.  Using Corollary 4.2, we
get
\begin{equation*} (x,y)= (\Delta_i(x), y(j)\ot y') =\sum_s(z_s, y(j))(z_s',y').\end{equation*}
  We are thus reduced to proving that 
$$(x,y) = 0 \ \  \mod (   q^{-1}\ba) \ \ \text{if}\  
  x\in {}_\ca\Lambda_i,\ \ y\in {}_\ca\Lambda_j,\ i\ne j.$$ It is
  obviously enough to prove this when $x$ and $y$ are products of the 
  $\tilde{P}_{k,i}$ and $\tilde{P}_{k,j}$, respectively.  Assume that $x
  = \tilde{P}_{k,i}$. If $y$ also has length one, this is just the
  statement of Proposition 3.1. Assume the result for all monomials $y$
  of length less than $s$. Write $y=\tilde{P}_{m,j}y_1$, where $y_1\in
 {}_\ca\Lambda_j$ has length less than $s$. We have, by Proposition
  \ref{comult}, that
\begin{equation*} (\tilde{P}_{k,i},\tilde{P}_{m,j}y_1)
  =(r(\tilde{P}_{k,i}), \tilde{P}_{m,j}\ot y_1)= \left(\sum_{l=0}^k
  \tilde{P}_{k-l,i}\ot \tilde{P}_{l,i}\ ,\ \tilde{P}_{m,j}\ot y_1\right)
\end{equation*} 
The right-hand side is zero mod $(q^{-1}\ba)$ by the induction hypothesis.
Assume now that we know the result for all monomials $y$ if $x$ has
length less than $s$. Write $x=x_1\tilde{P}_{k,i}$, where
$x_1\in{}_\ca\Lambda_i$ has length $s-1$.  Proceeding as before
and using Corollary 3.2, we see that
\begin{equation*}(x,y) = (x_1\ot\tilde{P}_{k,i}, r(y))=(x_1\ot\tilde{P}_{k,i}, \Delta_j(y)).\end{equation*}
The right-hand side is again zero mod $(q^{-1}\ba)$ by induction and the proof of (i) is complete.

To prove (ii), assume without loss of generality that $x(1)\ne 1$ and write $x=x(1)x'$.  Choose $j$ such that $y(j)\ne 1$ and $y(m)=1$ if $m<j$. If $j>1$ then, by Corollary 3.2, we see that 
$$r(y)\in\prod_{m\ne 1} {_\ca}\Lambda_m \ot {_\ca\Lambda}_m +{\text{ terms in $\bu^+(<)_+\bu^+(0)\ot\bu^+(0)\bu^+(>)_+$}}.$$ 
Hence, we get
$$(x,y) =(x(1)\ot x', r(y)) = 0\ \ \  \mod (   q^{-1}\ba)$$
by part (i) of Corollary \ref{orthogonal}. If $j=1$, write $y=y(1)y'$. If $y'=1$, then we are again done by part (i) of  Corollary \ref{orthogonal}. So assume that $y'(m)\ne 1 $ for some $m\ne 1$. Then, using Corollary 3.2, we can write
$$r(y) =y(1)\ot y' + 1\ot y(1)y'+\sum_sz_s\ot z_s'+ {\text{terms in }} \bu^+(<)\bu^+(0)\otimes
  \bu^+(0)\bu^+(>)_+,$$
where $z_s(m)\ne 1$. Since we have already proved that $(x(1), z_s(m))=0$
mod $(q^{-1}\ba)$ if $m\ne 1$, it follows that
$$(x,y)=(x(1),y(1))(x',y')\ \ \  \mod (   q^{-1}\ba).$$
The proof can now be completed by repeating the argument for $x'$ and $y'$. 
\hfill\qedsymbol


\section{Characterization of the canonical basis.}

The results of Section 3 allow us to use the theory symmetric functions
\cite{Ma} to modify the imaginary root vectors so that we have an
orthonormal basis mod $(q^{-1}\ba)$ for $\ua^+(0)$.  Using Theorem 2 and
Proposition \ref{luinner}, we then have an orthonormal basis mod $(q^{-1}\ba)$ for $\ua^+$. We use this basis to construct our crystal
basis  and the canonical basis.

Following \cite[page 41]{Ma}, we define, for $i\in I$ and a given a partition
$\lambda =(\lambda_1\ge\lambda_2\ge\cdots )$, the corresponding Schur
function $s_{\lambda,i}\in{}_\ca\Lambda_i$ by
\begin{equation*} s_{\lambda,i} = \det(\tilde{P}_{\lambda_k-k+m,i})_{1\le k,m\le t},\end{equation*}
  where $t\ge l(\lambda)$.  Next, given a function $\boc_0 \in
  \bn^{\car_0}$, consider the $n$-tuple of partitions $(\lambda^{(1)},
  \dots, \lambda^{(n)})$ whose $i$-th component is the partition with
  $\boc_0(k\delta,i)$ parts equal to $k$, and define $S_{\boc_0} =
  s_{{\lambda^{(1)}},1} s_{{\lambda^{(2)}},2} \dots
  s_{\lambda^{(n)},n}.$

\begin{defn} Let $\boc \in \bn^{\car}$.  Denote by
  $\boc_>$ (resp. $\boc_<$) the restriction of $\boc$ to $\car_>$ (resp.
  $\car_<$).  Define
\begin{equation} 
  B_\boc = (E_{\boc_>})\,.\, S_{\boc_0}\,.\, (E_{\boc_<}).
\end{equation}
\end{defn}

\begin{prop}  \label{biginner}
 \begin{enumerate}
\item[(i)]  The set $\{B_\boc\ | \ \boc\in \bn^{\car}\}$ is an $\A$-basis of $\ua^+$.
\item[(ii)] Let $\boc, \boc' \in \bn^{\car}$.  Then 
\begin{equation} \label{bloc} 
  (B_\boc, B_{\boc'})=\delta_{\boc,\boc'}\ \  
\mod (   q^{-1}\ba).
\end{equation}
\end{enumerate}
\end{prop}
\begin{pf} Let $\Lambda_i$ be the $\bz$-subring of $ {_\ca\Lambda_i}$
  generated by the $\tilde{P}_{k,i}$. The first statement is an
  immediate consequence of Theorem 2 and the fact \cite[Section 3.3]{Ma}
that the $s_{\lambda,i}$ form a $\bz$-basis of $ \Lambda_i$.   
  For the second, it suffices by Proposition \ref{luinner} to prove that
\begin{equation} \label{theend} (S_{\boc_0}, S_{\bf{c'_0}})=
  \delta_{\bf{c_0},\bf{c'_0}} \ \  
\mod (   q^{-1}\ba).
\end{equation} 
By Corollary
  \ref{orthogonal}, we see that
\begin{equation*}(S_{\boc_0}, S_{\bf{c'_0}})
  =(s_{{\lambda^{(1)}},1},s_{{\mu^{(1)}},1})
(s_{{\lambda^{(2)}},2},s_{{\mu^{(2)}},2})\dots
  (s_{{\lambda^{(n)}},n},s_{{\mu^{(n)}},n}) \ \  
\mod (   q^{-1}\ba).\end{equation*} Now we use \cite[Page 92, Exercise 25(c)]{Ma} and
  \cite[Chapter 1, Equation 4.8]{Ma} to conclude that
\begin{equation}
(s_{{\lambda^{(i)}},i},s_{{\mu^{(i)}},i}) 
= \delta_{\lambda^{(i)},\mu^{(i)}} \ \  
\mod (   q^{-1}\ba).\end{equation} 
  Part (ii) follows.
 \end{pf}

 Let $\bb$ be the canonical basis of $\bu^+$, and let $\cl$ be the
 $\bz[q^{-1}]$-lattice spanned by $\bb.$ By Proposition 5.1.3 in [K], we have
 the following alternative characterization of $\cl$:
\begin{equation} \cl = \{ x \in \ua^+\ |\ (x,x) \in \ba\}. \end{equation}

\begin{prop}  The set
\begin{equation} \{ B_\boc \ |\ \boc \in \bn^{\car} \}
\end{equation}
is a $\bz[q^{-1}]$-basis of $\cl.$
\end{prop}
\begin{pf}  This follows directly from Propositions \ref{luinner} and
\ref{biginner}, and   \cite[Lemma 16.2.5 (a)]{L2}.
\end{pf}

\begin{thm} \label{crystal} 
  For every $\boc \in \bn^{\car}$, there exists $b \in \bb$ such
  that
\begin{equation} \label{thm2} b = B_\boc  \ \  \mod (   q^{-1}\cl),
\end{equation}
and hence $\{B_\boc\ | \ \boc\in \bn^{\car}\}$ is a crystal basis of
$\ua^+$.
\end{thm}

We begin by observing that this theorem is certainly true up to sign:

\begin{prop} For every $\boc \in \bn^{\car}$, there exists $b \in \bb$ such
  that
\begin{equation}  b = \pm B_\boc  \ \  \mod (   q^{-1}\cl).
\end{equation}
\end{prop}
 \begin{pf} This follows from \cite[Lemma 16.2.5(f)]{L2} and Proposition \ref{biginner}.
\end{pf}

Let $\pi^0$ be the projection of $\ua^+$ onto $\ua^+(0)$ corresponding to the decomposition 
\begin{equation*} \ua^+ =\ua^+(0)\oplus (\ua^+(>)\ua^+(0)\ua^+(<)_+ +\ua^+(>)_+\ua^+(0)\ua^+(<)),\end{equation*}
which follows from Theorem 2. Note that this theorem also implies that the  restriction of $\pi^0$ to 
the $\A$-subalgebra of $\ua^+$ consisting of elements of homogeneity
in $\bn\delta$ is a homomorphism. The following additional properties of $\pi^0$ follow from Proposition 4.3:

\begin{cor} We have:
\begin{enumerate}
\item[(i)] $\cl$ is the $\bz[q^{-1}]$-lattice  spanned by $\{B_{\bold c}\}_{\bold c\in\bn^{\car}}$;
\item[(ii)] $\pi^0(\cl)\subseteq\cl$;
\item[(iii)] the non-zero elements of $\pi^0(\bold B)\subset\cl$ are linearly independent (mod $q^{-1}\cl$);
\item[(iv)] $\ua^+(0)\cap\cl$ is closed under products.
\end{enumerate}
\end{cor}
\begin{pf} Part (i) is immediate from Propositions 4.2 and 4.3. Part (ii) follows from part (i) and the fact that $\pi^0$ takes each $B_{\boc}$ either to itself or zero. Part (iii) follows from Proposition 4.3 and the argument in part (ii). Finally, by Proposition 4.2 and the definition of the $s_{\lambda,i}$, it follows that $\cl$ is the $\bz[q^{-1}]$-lattice spanned by $\{E_{\boc}\}_{\boc\in\bn^{\car}}$ (see Section 3). This implies that $\ua^+(0)\cap\cl$ is spanned over $\bz[q^{-1}]$ by the monomials in the $\tilde{P}_{k,i}$. The statement in (iv) is now clear.
\end{pf}

For the proof of Theorem 3, we note that by \cite[Proposition 8.3]{L3}, it suffices to prove the theorem for $\boc\in \bn^{\car_0}.$ This is done in the following lemmas. 
 
\begin{lem}  \label{lengthone} Let $i\in I$.  Then, for every $k > 0$,
\begin{equation*}
  \tilde P_{k,i}  =   \beta_{k,i} + b_{k,i},
\end{equation*}
where $b_{k,i} \in q^{-1}\cl$ and $\beta_{k,i} \in \bb.$
\end{lem}
\begin{pf} First we check that
\begin{equation} \label{firstmod}
  E_{\delta - \alpha_i}^{(k)} E_{\alpha_i}^{(k)} =   \tilde P_{k,i}
  \ \ \ \mod (   q^{-1}\cl).
\end{equation}
We know by Proposition \ref{greatidentity}
that 
\begin{equation} \label{keymod}
 E_{\delta - \alpha_i}^{(k)} E_{\alpha_i}^{(k)} = \tilde P_{k,i} + x,
\end{equation}
where $x \in \ua^+(>)\ua^+(0)\ua^+(<)_+.$
Since $r_i(E_{\delta - \alpha_i}^{(k)})= 0$, it follows that
\begin{equation*}
(E_{\delta - \alpha_i}^{(k)} E_{\alpha_i}^{(k)}\,, \,
E_{\delta - \alpha_i}^{(k)} E_{\alpha_i}^{(k)}) = 1 \ \ \mod (   q^{-1}\ba),
\end{equation*}
and considering \eqref{keymod}, this implies that
\begin{equation*}
      (\tilde{P}_{k,i},\tilde{P}_{k,i})+    (x,x) + 2(x,\tilde P_{k,i}) = 1 \ \ \mod (   q^{-1}\ba).
\end{equation*}
By Proposition 4.1, this means that 
\begin{equation*}
 (x,x) + 2(x,\tilde P_{k,i}) = 0 \ \ \mod (   q^{-1}\ba).
\end{equation*}
Proposition 4.3 implies that $(x,\tilde{P}_{k,i})=0$, and so finally 
$(x,x) = 0$ mod $(q^{-1}\ba)$. By \cite[Lemma 16.2.5(f)]{L2}, it
follows that
\begin{equation} \label{inlattice}  x \in q^{-1}\cl.
\end{equation}  
This proves \eqref{firstmod}.  Now, by
\cite[Proposition 8.2]{L3}, we have $E^{(k)}_{\delta - \alpha_i} = b$ mod
$(q^{-1} \cl$). Consider the Kashiwara operators $\tilde\phi_i: \ua^+\to\ua^+$ (introduced in \cite{Ka}) defined by
 \begin{equation*} \tilde\phi_i(E_{\alpha_i}^{(r)}x) = E_{\alpha_i}^{(r+1)}x, \end{equation*}  for all $r\ge 0$, $x\in\ua^+[i]$ (see Proposition \ref{ui}). Since $\sigma$ and $\tilde \phi_i$ map
$\cl$ into itself, and ${_ir}(E_{\delta-\alpha_i}) =0$ by \cite[Lemma 3.4]{B1},
we have
 \begin{equation} \label{kash} E_{\delta - \alpha_i}^{(k)}E_{\alpha_i}^{(k)} =
 \sigma \tilde \phi_i^k \sigma(E_{\delta - \alpha_i}^{(k)}) = \sigma
 \tilde \phi_i^k\sigma(b) \ \  \mod (   q^{-1} \cl),
\end{equation}
where $\sigma \tilde \phi_i^k \sigma(b) = b'$ mod $(q^{-1}\cl)$ and 
$b' \in \bb.$ The lemma now follows from \eqref{firstmod}.
\end{pf}

\begin{lem} \label{fixedi} For $i\in I$,  we have
\begin{equation} s_{\lambda,i} = \beta_{\lambda,i} + b_{\lambda,i},
\end{equation}
where $b_{\lambda,i} \in  q^{-1} \cl$ and
$\beta_{\lambda,i} \in \bb.$
\end{lem}
\begin{pf} We prove the lemma
  by induction on the length $\ell(\lambda)$ of $\lambda$.  The case 
  $\ell(\lambda) = 1$  is contained in  Lemma
  \ref{lengthone}.  Assume that the statement holds for all $\lambda$ such
  that $\ell(\lambda) < L$.  By the Pieri formulas \cite[Chapter 1]{Ma}, we have
\begin{equation} \label{pieri}
s_{(k),i} s_{\mu,i} = \sum_{\lambda \supset \mu} s_{\lambda,i},
\end{equation} 
where the summation is over those $\lambda$ such that $\lambda - \mu$ is
a horizontal $k$ strip.
By the inductive assumption, we have 
\begin{equation} s_{\mu,i} = \beta_{\mu,i} + b_{\mu,i}, \ s_{(k),i} = \beta_{(k),i} +
  b_{(k),i},
\end{equation} where $\beta_{\mu,i}, \beta_{(k),i} \in  \bb$ and  
$b_{\mu,i}, b_{(k),i} \in q^{-1} \cl$. Multiplying these expressions together and using Corollary 4.1, we obtain
\begin{equation} \label{firstreals}
s_{(k),i} s_{\mu,i} =  \beta_{(k),i} \beta_{\mu,i} + y,
\end{equation} 
where $\pi^0(y) \in q^{-1}\cl.$ By the positivity result for the 
canonical basis \cite[Theorem 14.4.13]{L2}, there exist $n_r \in
\bn[q,q^{-1}]$, $b_r \in \bb, \ r = 1, \dots, d$, such that
\begin{equation} \label{positivity1}
\beta_{(k),i} \beta_{\mu,i} = \sum_{r=1}^d n_r b_r.
\end{equation} 
It follows that 
\begin{equation} \label{comparison}
\sum_{\lambda \supset \mu} s_{\lambda,i} = \sum_{r=1}^d n_r b_r + y,
\end{equation}
where  $n_r \in \bn[q,q^{-1}]$. Hence,
\begin{equation}
\sum_{\lambda \supset \mu} s_{\lambda,i} = \sum_{r=1}^d n_r  \pi^{0}(b_r) +\pi^0(y).
\end{equation} 
 Since (4.7) holds up to sign, we have $s_{\lambda,i} = \pm \beta_{\lambda,i} +
b_{\lambda,i}$ with $b_{\lambda,i}\in q^{-1}\cl$, hence
\begin{equation*}
\sum_{\lambda \supset \mu} s_{\lambda,i} = \sum_{\lambda\supset\mu} \pi^{0}(\pm
\beta_{\lambda,i}) \ \ \ \ \  \mod (   q^{-1}\cl) 
\end{equation*}
and so
\begin{equation}
\sum_{\lambda\supset\mu} \pm\pi^0(\beta_{\lambda,i})=\sum n_r\pi^0(b_r)\ \ \ \ \ \mod (   q^{-1}\cl).
\end{equation}
On the right-hand side, we sum only over those $r$ for which $\pi^0(b_r)\ne 0$. We claim that, for such $r$, $n_r\in\bn[q^{-1}]$. For otherwise, we may assume that $q^N$ is the highest power of $q$ appearing in $n_{r_1},\ldots,n_{r_s}$, where $N,s>0$, and that the highest power of $q$ appearing in the other $n_r$ in (4.20) is $<N$. Multiplying (4.19) on both sides by $q^{-N}$ gives
\begin{equation*}
\sum_{t=1}^s({q^{-N}n_{r_t}})_\infty\pi^0(b_{r_t})=0\ \ \ \ \ \mod (   q^{-1}\cl),
\end{equation*}
where $(\zeta)_\infty$ denotes the constant coefficient of an element $\zeta\in\bn[q^{-1}]$. This contradicts the linear independence of the $\pi^0(b_{r_t})$. Thus,  
\begin{equation*} \sum_{\lambda\supset\mu}\pm\pi^{0}(\beta_{\lambda,i})=\sum_r({n_r})_\infty\pi^0(b_r).\end{equation*}
It follows from linear independence again that, for  each $r$ such that $\pi^0(b_r)\ne 0$, there exists a $\lambda\supset\mu$ such that $\beta_{\lambda,i} =b_r$ and $({n_r})_\infty=\pm 1$. Hence, all the signs must be $+$.

Since every $\lambda$ of length $L$ appears in an equation of the form
\eqref{pieri}, this implies the lemma.
\end{pf}

\begin{lem} \label{differenti} 
Given an ``imaginary'' PBW basis monomial $S_{\boc_0}$, we have 
\begin{equation}
S_{\boc_0} = \beta_{\boc_0} + b_{\boc_0},
\end{equation}
where $b_{\boc_0} \in q^{-1} \cl$ and
$\beta_{\boc_0} \in \bb.$
\end{lem}
\begin{pf} Write
\begin{equation}
S_{\boc_0} = \prod_i s_{\lambda^{(i)},i} = \pm \beta_{\boc_0} + y,
\end{equation}
where $\beta_{\boc_0} \in \bb$ and $y \in q^{-1}\cl.$  Also
write each 
\begin{equation}
s_{\lambda^{(i)},i}  = \beta_{\lambda^{(i)}} + b_{\lambda^{(i)}},
\end{equation}
where by Lemma \ref{fixedi}, $\beta_{\lambda^{(i)}} \in \bb$ and 
$b_{\lambda^{(i)}} \in  q^{-1} \cl.$
By the positivity result for the canonical basis \cite[Theorem 14.4.13]{L2}
and (i)--(iv) above, we have 
\begin{equation}
\prod_i s_{\lambda^{(i)},i} = \sum_{r=1}^k n_r \beta_r + y,
\end{equation}
where $n_r \in \bn[q,q^{-1}]$ and $\pi^0(y) \in q^{-1}\cl .$
Arguing exactly  as in Lemma \ref{fixedi}, the statement follows.
\end{pf}

To obtain the canonical basis, let $\pi: \cl \rightarrow \cl/q^{-1}\cl$
be the natural projection. Then, $\pi$ takes the basis $\{B_\boc\ | \ 
\boc\in \bn^{\car}\}$ to a $\bz$-basis of $\cl / q^{-1} \cl$ and it
is known \cite{Ka} that $\pi$ restricts to an isomorphism of
$\Z$-modules $\pi': \cl \cap \ov{\cl} \isom \cl / q^{-1} \cl.$
We have the following immediate corollary of Theorem \ref{crystal}.

\begin{thm}      
\begin{equation*} 
  \bb = \{{\pi'}^{-1}\pi(B_\boc) \ | \ \boc \in \bn^{\car} \} .\ \ \ \ \ \ \ \ \ \ \ \ \ \ \qedsymbol
\end{equation*}
\end{thm}

\end{document}